\documentclass[11pt]{amsart}
\usepackage{amsmath}
\usepackage{amsthm}
\usepackage{amscd}
\usepackage{amssymb}
\usepackage{amsfonts}
\usepackage{graphicx}
\newtheorem{theorem}{Theorem}[section]
\newtheorem{lemma}[theorem]{Lemma}
\newtheorem{corollary}[theorem]{Corollary}
\newtheorem{proposition}[theorem]{Proposition}

\newtheorem{definition}[theorem]{Definition}
\theoremstyle{definition}
\newtheorem{example}[theorem]{Example}

\theoremstyle{remark}
\newtheorem*{remark}{Remark}

\newtheorem*{ack}{Acknowledgments}

\numberwithin{equation}
{section}
\newcounter{temp}

\makeatletter
\def\square{\RIfM@\bgroup\else$\bgroup\aftergroup$\fi
\vcenter{\hrule\hbox{\vrule\@height.6em\kern.6em\vrule}\hrule}\egroup}
\makeatother

\DeclareMathOperator{\Specf}{Specf}
\DeclareMathOperator{\Spec}{Spec}
\DeclareMathOperator{\Proj}{Proj}
\DeclareMathOperator{\Hom}{Hom}
\DeclareMathOperator{\Der}{Der}

\DeclareMathOperator{\Ker}{Ker}
\DeclareMathOperator{\Coker}{Coker}
\DeclareMathOperator{\cone}{cone}
\DeclareMathOperator{\st}{st}
\DeclareMathOperator{\link}{lk}
\DeclareMathOperator{\flip}{Fl}

\DeclareMathOperator{\Ann}{Ann}
\DeclareMathOperator{\Def}{Def}

\DeclareMathOperator{\Ext}{Ext}

\DeclareMathOperator{\id}{id}
\DeclareMathOperator{\Sym}{Sym}

\begin{document}

\title{Deforming Stanley-Reisner schemes}
\author{Klaus Altmann \and Jan Arthur Christophersen}
\maketitle
\begin{abstract} We study the deformation theory of projective
Stanley-Reisner schemes associated to combinatorial manifolds. We achieve
detailed descriptions of first order deformations and obstruction spaces.
Versal base spaces are given for certain Stanley-Reisner surfaces.
\end{abstract}
    
\section{Introduction} 
We consider the deformation theory of projective Stanley-Reisner schemes
associated to combinatorial manifolds. This paper builds on the results of
\cite{ac:cot} where we described the cotangent cohomology of
Stanley-Reisner rings for arbitrary simplicial complexes. 

Smoothings of Stanley-Reisner schemes associated to combinatorial manifolds
yield interesting algebraic geometric varieties. For example if the complex
is a triangulated sphere then the smoothing (if possible) would be
Calabi-Yau. The Stanley-Reisner scheme of a triangulated torus would smooth
to an abelian variety. A triangulated $\mathbb{RP}^{2}$ would give an
Enriques surface. It is our hope that the results of this paper may be
useful for the study of degenerations of such special varieties.

In the surface case there will be non-algebraic deformations of these
Stanley-Reisner schemes. To separate the algebraic deformations we use the
functor $\Def_{(X,L)}$ of deformations of the pair $(X,L)$, $X$ a scheme
and $L$ an invertible sheaf on $X$. In Section~\ref{DXL} we state and prove
properties of this functor for singular schemes.

We can give a very explicit account of first order deformations and
obstruction spaces. In the curve, surface and threefold case we are able to
give dimension formulas. This is done in Sections~\ref{tiA} and
\ref{c_def}. 

In the surface case we detail the non-algebraic deformations in the
beginning of Section~\ref{algdef}. We conclude the paper with a description
of the versal base space of algebraic deformations for $2$-dimensional
combinatorial manifolds with vertex valencies not greater than $6$.

\begin{ack} We are grateful to Edoardo Sernesi for several discussions, in
particular related to the functor $\Def_{(X,L)}$. We would also like to
thank Jan Stevens for pointing out an error in the first version. The paper
was finished while the second author was on sabbatical visiting Freie
Universit\"{a}t Berlin where he received support for his stay.
\end{ack}
\section{Preliminaries}

\subsection{Simplicial complexes and combinatorial manifolds} \label{simp}
Let $[n]$ be the set $\{0,\ldots,n\}$ and let $\Delta_{n}:=2^{[n]}$ be the
full simplex. 
A {\em simplicial complex} for us is a subset $\mathcal{K}\subseteq
\Delta_{n}$ satisfying the face relation: $f \in \mathcal{K}\,  \& \,  g
\subseteq f
\Rightarrow g \in \mathcal{K}$. We denote the the support of $\mathcal{K}$
by 
$[\mathcal{K}] = \{i \in [n] \, | \, \{i\} \in \mathcal{K}\}$.

For $g \subseteq [n]$, denote by $\bar{g}:=2^g$ and $
\partial g := \bar{g}\setminus \{g\}$ the full simplex and its boundary,
respectively. The {\it join} $\mathcal{K}\ast \mathcal{L}$ of two complexes
$\mathcal{K}$ and $\mathcal{L}$ is the complex defined by
\vspace{-1ex}
\[
\mathcal{K}\ast \mathcal{L} := \{f\vee g : f\in \mathcal{K},\, g\in
\mathcal{L}\}
\]
where $\vee$ means the disjoint union. If $f\in \mathcal{K}$ is a face, we
may define
\begin{itemize}
\item the {\it link} of $f$ in $\mathcal{K}$;
$\;\link(f,\mathcal{K}):=\{g\in \mathcal{K}: g\cap f =
\emptyset \text{ and } g\cup f\in \mathcal{K}\}$,

\item the {\it open star} of $f$ in $\mathcal{K}$; 
$\;\st(f,\mathcal{K}):=\{g\in \mathcal{K}: f\subseteq g\}$, and

\item the {\it closed star} of $f$ in $\mathcal{K}$; 
$\;\overline{\st}(f,\mathcal{K}):=\{g\in \mathcal{K}: g\cup f\in
\mathcal{K}\}$.

\end{itemize}
Notice that the closed star is the subcomplex
$\overline{\st}(f,\mathcal{K}) =
\bar{f}\ast \link(f,\mathcal{K})$.
% \end{nota}
%
The {\em geometric realization} of $\mathcal{K}$, denoted $|\mathcal{K}|$,
is defined as
\[
|\mathcal{K}| = \big\{\alpha: [n] \to [0,1] | \{i | \alpha(i) \ne 0\}\in
\mathcal{K} 
\mbox{ and  $\,\sum_i \alpha(i) = 1$} \big\}\, .
\]
To every non-empty $f\in \mathcal{K}$, one assigns the {\em relatively
open} simplex $\langle f\rangle \subseteq |\mathcal{K}|$;
\[
\langle f\rangle  = \{\alpha\in |\mathcal{K}| \, | \, \alpha(i) \ne 0
\text{ if and only if } i\in f \}\, .
\]
On the other hand, each subset $Y \subseteq \mathcal{K}$, i.e. $Y$ is not
necessarily a subcomplex, determines a topological space
\[
\langle Y\rangle:=
\begin{cases}
        \bigcup_{f\in Y}\langle f\rangle& \text{if $\emptyset\not\in Y$},
\\
        \cone \left(\bigcup_{f\in Y}\langle f\rangle\right)& \text{if
        $\emptyset\in Y$}\, .
\end{cases}
\]
In particular, $\langle \mathcal{K}\setminus\{\emptyset\}\rangle =
|\mathcal{K}|$ and $\langle
\mathcal{K}\rangle = |\cone(\mathcal{K})|$ where $\cone(\mathcal{K})$ is
the simplicial complex $\Delta_0\ast \mathcal{K}$.

If $f$ is an $r$-dimensional face of $\mathcal{K}$, define the {\em
valency} of $f$, $\nu(f)$, to be the number of $(r+1)$-dimensional faces
containing $f$.  Thus $\nu(f)$ equals the number of vertices in
$\link(f,\mathcal{K})$.

In this paper we are mostly interested in combinatorial manifolds. We refer
to  \cite{hud:pie} for definitions and results in $PL$ topology.  A {\em
combinatorial $n$-sphere} is a simplicial complex $\mathcal{K}$ such that
$|\mathcal{K}|$ is $PL$-homeomorphic to $|\partial \Delta_{n+1}|$.  A
simplicial complex $\mathcal{K}$ is a {\em combinatorial $n$-manifold} if
for all non-empty faces $f\in \mathcal{K}$, $|\link(f,\mathcal{K})|$ is a
combinatorial sphere of dimension $n-\dim f -1$.  If we also allow
$|\link(f,\mathcal{K})|$ to be a ball of dimension $n-\dim f -1$, then
$\mathcal{K}$ is called a combinatorial manifold with boundary. In this
case we denote the boundary $\partial \mathcal{K} = \{f \in \mathcal{K}
\, | \, |\link(f,\mathcal{K})| \text{ is a ball} \}$. In dimensions less
than four all triangulations of topological manifolds are combinatorial
manifolds (see e.g. \cite{hud:pie}). In this paper we call $\mathcal{K}$ a
{\em manifold} if it is a combinatorial manifold without boundary.

We will need notation for some special manifolds. Write $\Sigma
\mathcal{K}$ for the suspension of a complex $\mathcal{K}$.  Let $E_n$ be
the boundary of the $n$-gon; i.e.\ $|E_n|\approx S^1$. Let $C_n$ be the
chain of $n$ $1$-simplices; i.e.\ $|C_n|\approx B^1$.  Let $\partial C(n,3)
=
\partial\bigtriangleup_1 \ast  C_{n-3}
\cup
\partial  C_{n-3} \ast \bigtriangleup_1$ be the boundary of the
$3$-dimensional cyclic polytope (see \cite[4.7]{gr:con}).  If
$[\bigtriangleup_1]=\{0,n-1\}$ and $[C_{n-3}]=\{1, 2,
\dots , n-2\}$ then the facets of $
\partial C(n,3)$ are
\begin{multline*}
\{0,2,n-1\}, \{0,n-2,n-1\}, \{0,2,3\}, \{0,3,4\}, \dots , \{0,n-3,n-2\},\\
\{2,3,n-1\}, \{3,4,n-1\}, \dots , \{n-3,n-2,n-1\}.
\end{multline*} A drawing of this complex for $n=7$ may be found in
Section~\ref{c_def}.

\subsection{Stanley-Reisner schemes}\label{sr_rings} Let
$P=k[x_0,\ldots,x_n]$ be the polynomial ring in $n+1$ variables over an
algebraically closed field field $k$.  If $a=\{i_{1},\dots ,i_{k}\}\in
\Delta_{n}$, we write $x_a\in P$ for the square free monomial
$x_{i_1}\cdots x_{i_k}$.  If $\mathbf{a}=(a_0,\dots, a_n)\in
\mathbb{Z}^{n+1}$, set $x^\mathbf{a}\in P$ to be the monomial
$x_0^{a_{0}}\cdots x_n^{a_n}$. The support of $\mathbf{a}$ is defined as
$a:=\{i\in[n] \, | \,   \mathbf{a}_i\neq 0\}$. We will throughout write
 $\mathbf{c}=\mathbf{a}-\mathbf{b}$ for the decomposition of $\mathbf{c}$
in its positive and negative part, i.e.,
$\,\mathbf{a},\mathbf{b}\in\mathbb{N}^{n+1}$ with both elements having
disjoint supports $a$ and $b$, respectively.

A simplicial complex $\mathcal{K}\subseteq
\Delta_{n}$ gives rise to an ideal
\[
I_\mathcal{K}:=\langle x_{p} \, | \,  p\in \Delta_{n}\setminus
\mathcal{K}\rangle
\subseteq P .
\]
The {\it Stanley-Reisner ring} is then $A_\mathcal{K}=P/I_\mathcal{K}$.  We
refer to \cite{sta:com} for more on Stanley-Reisner rings. 

We can associate the schemes ${\mathbb A}(\mathcal{K})=\Spec A_\mathcal{K}$
and ${\mathbb P}(\mathcal{K}) = \Proj A_\mathcal{K}$ with these rings. The
latter looks like $\lvert \mathcal{K} \rvert$ -- its simplices have just
been replaced by projective spaces. If $f$ is a subset of $[n]$, let
$D_{+}(x_f)
\subseteq {\mathbb P}(\mathcal{K})$ be the chart corresponding to
homogeneous localization of $A_\mathcal{K}$ by the powers of $x_f$. Then
$D_{+}(x_f)$ is empty unless $f\in \mathcal{K}$ and if $f\in \mathcal{K}$
then $$
D_{+}(x_f) = {\mathbb A}(\link(f,\mathcal{K})) \times (k^{*})^{\dim f} \,
.$$

We will need the following result of Hochster as stated in \cite[Proof of
Theorem 4.1]{sta:com}.
\begin{theorem}\label{loco} Let $\mathfrak{m}$ be the irrelevant maximal
ideal in the multi-graded ring $k[x_0,\dots x_n]$.  Let
$H^i_{\mathfrak{m}}(A_\mathcal{K})_{\mathbf{c}}$ be a multi-graded piece of
the local cohomology module with $\mathbf{c} \in \mathbb{Z}^n$.  Then
$H^i_{\mathfrak{m}}(A_\mathcal{K})_{\mathbf{c}}=0$ unless $\mathbf{c}\le
\mathbf{0}$, i.e.\ $\mathbf{c} = \mathbf{0} - \mathbf{b}$, and $b \in
\mathcal{K}$ in which case $$
H^i_{\mathfrak{m}}(A_\mathcal{K})_{\mathbf{c}} \simeq
\widetilde{H}^{i-|b|-1}(\link(b);k) \, .$$
\end{theorem} Recall that by comparing the \v{C}ech complex of
$\bigoplus_{m}
\mathcal{O}_{\Proj A}(m)$ and the complex computing $H^i_{\mathfrak{m}}(A)$
we get $\bigoplus_{m} H^i(\Proj A,
\mathcal{O}_{\Proj A}(m)) \simeq H^{i+1}_{\mathfrak{m}}(A)$ when $i\ge 1$
and an exact sequence $$
0 \to H^{0}_{\mathfrak{m}}(A) \to A \to \bigoplus_{m} H^0(\Proj A,
\mathcal{O}_{\Proj A}(m)) \to H^{1}_{\mathfrak{m}}(A) \to 0 \, .$$
As a consequence we get
\begin{theorem}\label{O_cohom} If $\mathcal{K}$ is a simplicial complex 
then $$
H^p({\mathbb P}(\mathcal{K}),\mathcal{O}_{{\mathbb P}(\mathcal{K})})
\simeq H^{p}(\mathcal{K};k)$$
and if $m\ge 1$ 
$$
H^p({\mathbb P}(\mathcal{K}), \mathcal{O}_{{\mathbb P}(\mathcal{K})}(m)) =
\begin{cases}  (A_{\mathcal{K}})_{m} & \text{if $p = 0$ }\\
0 & \text{if $p \ge 1$ }
\end{cases}$$
\end{theorem}

\subsection{The cotangent spaces and sheaves}\label{cot} For standard
definitions and results in deformation theory of schemes we refer to
\cite{ser:def}. To fix notation we recall that for an $S$-algebra $A$ and
an $A$-module $M$ there exist the cotangent modules $T^{i}_{A/S}(M)$. We
write $T^{i}_{A}$ when  $S=k$ and $M=A$. The module $T^{0}_{A} =
\Der_{k}(A,A)$ consists of the infinitesimal automorphisms of $A$,
$T^{1}_{A} \simeq \Def_{\Spec A}(k[\epsilon])$ is the space of first order
deformations of $\Spec A$ and $T^{2}_{A}$ contains the obstructions for
lifting deformations.

If $Y$ is a scheme we may globalize these modules. (See for example
\cite[Appendice]{an:hom} and \cite[3.2]{la:for}.) Let $\mathcal S$ be a
sheaf of rings on $Y$, $\mathcal A$ an $\mathcal S$ algebra and $\mathcal
F$ an ${\mathcal A}$ module.  We get the cotangent cohomology sheaves
${\mathcal T}_{{\mathcal A}/{\mathcal S}}^i({\mathcal F})$ as the sheaves
associated to the presheaves $U\mapsto T^i({\mathcal A}(U)/{\mathcal
S}(U);{\mathcal F}(U))$.

There are also the groups $T_{{\mathcal A}/{\mathcal S}}^i({\mathcal F})$ -
the hyper-cohomology of the cotangent complex on $Y$.  If ${\mathcal
A}={\mathcal F} = {\mathcal O}_Y$ and $S=k$, then (abbreviating as above)
the $T^i_Y$ play the same role in the deformation theory of $Y$ as in the
local case. There is a ``local-global" spectral sequence $$
E^{p,q}_{2} = H^p(Y,{\mathcal T}_Y^q)\Rightarrow T_Y^{p+q}$$
which relates the local and global deformations.  In particular first order
automorphisms are described as $T_Y^0 = H^0(Y,\Theta_Y)$ and there is an
exact sequence $$
0 \to H^1(Y,{\mathcal T}_Y^0) \to T_Y^1 \to H^0(Y,{\mathcal T}_Y^1)
\to H^2(Y,{\mathcal T}_Y^0)\, .$$
All three groups $H^0(Y,{\mathcal T}_Y^2)$, $H^1(Y,{\mathcal T}_Y^1)$ and
$H^2(Y,{\mathcal T}_Y^0)$ contribute to the obstructions.

\section{The functor $\Def_{(X,L)}$}\label{DXL}

Let $X$ be a scheme over an algebraically closed field $k$ and $L$ an
invertible sheaf on $X$. Let $A$ be an object in the category $\mathcal{A}$
of local artinian $k$-algebras with residue field $k$. We recall the
definition of the functor $\Def_{(X,L)}$ of infinitesimal deformations of
the pair $(X,L)$ in \cite[3.3.3]{ser:def} and generalize its properties to
singular schemes.

An infinitesimal deformation of the pair $(X,L)$ over $A$ is a deformation
$\mathcal{X} \to \Spec(A)$ with an invertible sheaf $\mathcal{L}$ on
$\mathcal{X}$ such that $\mathcal{L}_{|X} = L$. Two such deformations
$(\mathcal{X}, \mathcal{L})$ and $(\mathcal{X}^{\prime}, 
\mathcal{L}^{\prime})$ are isomorphic if there is an isomorphism of
deformations $f : \mathcal{X} \to \mathcal{X}^{\prime}$ and an isomorphism
$\mathcal{L} \to f^{*}\mathcal{L}^{\prime}$. Let $\Def_{(X,L)} :
\mathcal{A}\to (\text{sets})$ denote the corresponding functor of Artin
rings. We define $\Def^{\prime}_{(X,L)}$ to be the subfunctor of
deformations of the pair where the deformation of $X$ is {\em locally
trivial}.

For any scheme there is a natural map $\mathcal{O}_{X}^{*} \to
\Omega^{1}_{X}$ defined locally by 
$$
u \mapsto \frac{du}{u} \, .$$
Let $c: H^{1}(X, \mathcal{O}_{X}^{*}) \to H^{1}(X, \Omega^{1}_{X})$ be the
induced map in cohomology. Now $ H^{1}(X, \Omega^{1}_{X}) \simeq
\Ext^{1}(\mathcal{O}_{X}, \Omega^{1}_{X})$, so $c(L)$ gives us an extension
\begin{equation*}
\label{ex} e_{L}: \quad 0 \to \Omega^{1}_{X} \to \mathcal{Q}_{L} \to
\mathcal{O}_{X} \to 0 \, . 
\end{equation*} In the smooth case $\mathcal{P}_{L} = \mathcal{Q}_{L}
\otimes_{\mathcal{O}_{X}} L$ is known as the sheaf of principle parts of
$L$.

Set $\mathcal{E}_{L} := \mathcal{Q}_{L}^{\vee}$ and note that the dual
sequence 
$$
0 \to \mathcal{O}_{X} \to \mathcal{E}_{L} \to \Theta_{X} \to 0 $$
is also exact. In the smooth case this is known as the Atiyah extension
associated to $L$.

We generalize \cite[Theorem 3.3.11]{ser:def}.
\begin{theorem} \label{DXL_theorem} Let $X$ be a reduced projective scheme
and $L$ an invertible sheaf on $X$. Then:
\begin{list}{\textup{(\roman{temp})}}{\usecounter{temp}}
\item The functor $\Def_{(X,L)}$ has a hull.
\item There are isomorphisms $\Def_{(X,L)}(k[\epsilon]) \simeq
\Ext^{1}_{\mathcal{O}_{X}}( \mathcal{Q}_{L}, \mathcal{O}_{X})$ and
$\Def^{\prime}_{(X,L)}(k[\epsilon]) \simeq H^1(X,\mathcal{E}_{L})$ and an
exact sequence of $k$-vector spaces 
$$
0 \to H^1(X, \mathcal{E}_{L}) \to \Ext^{1}_{\mathcal{O}_{X}}(
\mathcal{Q}_{L},
\mathcal{O}_{X}) \to H^0(X, \mathcal{T}^{1}_{X}) \to H^2(X,
\mathcal{E}_{L}) \, .$$
\item The obstructions for $\Def_{(X,L)}$ lie in $H^0(X,
\mathcal{T}^{2}_{X})$, $H^{1}(X, \mathcal{T}^{1}_{X})$ and $H^2(X,
\mathcal{E}_{L})$.
\item Given a first-order deformation of $X$ with isomorphism class $\xi
\in \Ext^{1}(\Omega^{1}_{X}, \mathcal{O}_{X})$, there is a first-order
deformation of $L$ along $\xi$ if and only if in the Yoneda product 
$$
\Ext^{1}(\Omega^{1}_{X}, \mathcal{O}_{X}) \times
\Ext^{1}(\mathcal{O}_{X}, \Omega^{1}_{X}) \to \Ext^{2}(\mathcal{O}_{X},
\mathcal{O}_{X}) = H^{2}(X, \mathcal{O}_{X})$$
we have $\xi\cdot c(L) = 0$.
\item If $L$ is very ample and $H^{1}(X,L) = 0$ then any formal deformation
of the pair $(X,L)$ is effective.
\end{list}
\end{theorem}

\begin{remark} It follows from (i) and (v) and a theorem of Artin
(\cite[Theorem 2.5.14]{ser:def}) that under the conditions in (v),
$\Def_{(X,L)}$ has an {\em algebraic} versal deformation.
\end{remark}

\begin{proof} In the proof of \cite[Theorem 3.3.11]{ser:def} the
Schlessinger conditions are checked for $\Def_{(X,L)}$ in the case $X$ is
nonsingular, but nowhere is the assumption nonsingular needed. 

For the remainder of the proof  choose an affine cover $\{U_{i}\}$ of $X$.
Let $L$ be  represented by a \v{C}ech cocycle $(f_{ij})$, $f_{ij} \in
\Gamma(U_{ij},
\mathcal{O}^{*}_{X})$. 

(ii) We will define a map $\Phi :\Def_{(X,L)}(k[\epsilon]) \to
\Ext^{1}_{\mathcal{O}_{X}}( \mathcal{Q}_{L}, \mathcal{O}_{X})$. Recall
first the isomorphism  $\Def_{X}(k[\epsilon]) \to
\Ext^{1}_{\mathcal{O}_{X}}(
\Omega^{1}_{X}, \mathcal{O}_{X})$  in the reduced case. If $\mathcal{X} \to
\Spec(k[\epsilon])$ is a first-order deformation, then the cotangent
sequence for $k \to \mathcal{O}_{\mathcal{X}} \to \mathcal{O}_{X}$ becomes
the exact sequence 
$$
0 \to \mathcal{O}_{X} \to \Omega^{1}_{\mathcal{X}} \otimes_{k[\epsilon]} k
\to \Omega^{1}_{X} \to 0$$
and the class of this extension in
$\Ext^{1}_{\mathcal{O}_{X}}(\Omega^{1}_{X}, \mathcal{O}_{X})$ is the image
of the isomorphism class of $\mathcal{X}$.

If $(\mathcal{X}, \mathcal{L})$ represents a first-order deformation we may
construct an extension $e_{\mathcal{L}}$: $$
0 \to \Omega^{1}_{\mathcal{X}} \to \mathcal{Q}_{\mathcal{L}} \to
\mathcal{O}_{\mathcal{X}} \to 0 $$
and a commutative diagram of exact sequences 
$$
\begin{CD} 0 @>>> \Omega^{1}_{\mathcal{X}} \otimes_{k[\epsilon]} k @>>>
\mathcal{Q}_{\mathcal{L}} \otimes_{k[\epsilon]} k @>>> \mathcal{O}_{X} @>>>
0\\
@. @VV\alpha V @VV\beta V @VV=V @.\\
0 @>>> \Omega^{1}_{X} @>>> \mathcal{Q}_{L} @>>> \mathcal{O}_{X} @>>> 0
\end{CD}$$
with surjective vertical maps. Thus  $\ker(\beta) \simeq \ker(\alpha)
\simeq \mathcal{O}_{X}$. This yields an exact sequence 
$$
0 \to \mathcal{O}_{X}  \to  \mathcal{Q}_{\mathcal{L}} \otimes_{k[\epsilon]}
k \to \mathcal{Q}_{L}
\to 0 $$
defining $\Phi$.

To describe $\Phi^{-1}$ we look again at why 
$\Ext^{1}_{\mathcal{O}_{X}}(\Omega^{1}_{X}, \mathcal{O}_{X}) \simeq
\Def_{X}(k[\epsilon])$. If 
$$
0 \to \mathcal{O}_{X}  \to  \mathcal{A}\stackrel{p}{\to} \Omega^1_{X} \to
0$$
defines an element of $\Ext^{1}_{\mathcal{O}_{X}}(\Omega^{1}_{X},
\mathcal{O}_{X})$, then construct the first-order deformation with
structure sheaf $\mathcal{O}_{\mathcal{X}}:= \mathcal{A}
\times_{\Omega^1_{X}} \mathcal{O}_{X}$, where the fibre product is with
respect to $p$ and the universal derivation $d : \mathcal{O}_{X} \to
\Omega^1_{X}$. One can then show that $\mathcal{A} \simeq
\Omega^{1}_{\mathcal{X}} \otimes_{k[\epsilon]} k$.

Over an open $U \subset X$, $\mathcal{O}_{\mathcal{X}}$ is the
$k[\epsilon]$ algebra $\{f +
\epsilon  a : (a,f) \in \Gamma(U, \mathcal{A} \times_{\Omega^1_{X}}
\mathcal{O}_{X})\}$. Note that the units $\Gamma(U,
\mathcal{O}^{*}_{\mathcal{X}}) = \{f + \epsilon a \in \Gamma(U,
\mathcal{O}_{\mathcal{X}}) : f \in \Gamma(U,
\mathcal{O}^{*}_{X})\}$.

Now let $$
0 \to \mathcal{O}_{X}  \to  \mathcal{B} \stackrel{q}{\to} \mathcal{Q}_{L}
\to 0$$
define an element of $\Ext^{1}_{\mathcal{O}_{X}}(\mathcal{Q}_{L},
\mathcal{O}_{X})$. From the extension $e_{L}$ we have a map $\alpha:
\Omega^1_{X} \to  \mathcal{Q}_{L}$ and we may construct the pullback
extension by $\alpha$. Let the middle term in this extension be
$\mathcal{A} = \mathcal{B} \times_{ \mathcal{Q}_{L}} \Omega^1_{X}$. We get
a commutative diagram with exact rows and columns:
\begin{equation}
\label{B}
\begin{CD} 
@.  @. 0 @. 0 @.\\
@. @. @VVV @VVV @.\\
0 @>>> \mathcal{O}_{X} @>>> \mathcal{A} @>p>> \Omega^{1}_{X} @>>> 0\\
@. @VV=V @VVV @VV\alpha V @.\\
0 @>>>  \mathcal{O}_{X} @>>> \mathcal{B} @>q>> \mathcal{Q}_{L} @>>> 0\\
@. @. @VVV @VVV @.\\
@.  @.\mathcal{O}_{X} @>=>> \mathcal{O}_{X} @.\\
@. @. @VVV @VVV @.\\
@.  @. 0 @. 0 @.
\end{CD}
\end{equation}
where the right column is $e_{L}$ and the the first row defines a first
order deformation $\mathcal{O}_{\mathcal{X}}$ as above .

To create $\mathcal{L}$ we need a cocycle  $(F_{ij})$, $F_{ij} \in
\Gamma(U_{ij},
\mathcal{O}^{*}_{\mathcal{X}})$ lifting the $(f_{ij})$. That means $F_{ij}
= f_{ij} + \epsilon a_{ij}$, $a_{ij} \in \Gamma(U_{ij}, \mathcal{A})$ with
$p(a_{ij}) = d f_{ij}$. The cocycle condition $F_{ij}F_{jk}=F_{ik}$ may be
computed to be equivalent to $$
\frac{a_{ij}}{f_{ij}} + \frac{a_{jk}}{f_{jk}} = \frac{a_{ik}}{f_{ik}}\, .$$
Thus $b_{ij} = a_{ij}/f_{ij}$ defines a class in $H^{1}(X, \mathcal{A})$
and 
$$
p(b_{ij}) = \frac{df_{ij}}{f_{ij}} = [e_{L}] \in H^{1}(X, \Omega^1_{X})\,
.$$

So to construct $\mathcal{L}$ we need to find a class in $p^{-1}(e_{L})
\subseteq H^{1}(X, \mathcal{A})$. A diagram chase shows that $e_{L}$ is the
pushout of the middle column of the diagram
\ref{B} by $p$. Thus the extension class of 
$$
0 \to \mathcal{A}  \to  \mathcal{B} \to \mathcal{O}_{X} \to 0$$
in $H^{1}(X, \mathcal{A})$ give us the wanted class. To be precise this
class is $\delta(1)$ where $\delta : H^0(\mathcal{O}_{X}) \to
H^{1}(\mathcal{A})$ is induced from the exact sequence. This also shows
that this extension is $e_{\mathcal{L}} \otimes k$ so we have defined
$\Phi^{-1}$.

The local-global spectral sequence for $\Ext$ yields a  four-term exact
sequence \begin{multline*} 0 \to H^1(X, \mathcal{E}_{L}) \to
\Ext^{1}_{\mathcal{O}_{X}}(
\mathcal{Q}_{L},
\mathcal{O}_{X})  \\
\to H^0(X, Ext^{1}_{\mathcal{O}_{X}}(\mathcal{Q}_{L},
\mathcal{O}_{X}))
\to H^2(X, \mathcal{E}_{L}) 
\end{multline*} which is almost what we want. Apply
$Ext(-,\mathcal{O}_{X})$ to $e_{L}$ to get $$
Ext^{1}_{\mathcal{O}_{X}}(\mathcal{Q}_{L},
\mathcal{O}_{X}) \simeq Ext^{1}_{\mathcal{O}_{X}}(\Omega^{1}_{X},
\mathcal{O}_{X})  \simeq \mathcal{T}^{1}_{X} \, .$$ 
This proves the existence of the exact sequence in (ii). 

(iii) Consider a small extension 
$$
0 \to (t) \to A^{\prime} \to A \to 0$$
of local artinian $k$-algebras and let $(\mathcal{X}, \mathcal{L})$ be a
deformation over $A$. The obstructions in the two first spaces are well
known. If they vanish we are in the following situation:
\begin{list}{\textup{(\alph{temp})}}{\usecounter{temp}}
\item On each $U_{i}$ we have deformations $(U_{i},
\mathcal{O}^{\prime}_{i}) \to \Spec(A^{\prime})$ of the affine schemes
$(U_{i}, \mathcal{O}_{X}|_{U_{i}})$ lifting $(U_{i},
\mathcal{O}_{\mathcal{X}}|_{U_{i}})$.
\item  On each $U_{ij}$ we have isomorphisms $\phi_{ij} :
\mathcal{O}^{\prime}_{i}|_{U_{ij}} \to
\mathcal{O}^{\prime}_{j}|_{U_{ij}}$ lifting the identity on
$\mathcal{O}_{\mathcal{X}}|_{U_{ij}}$. Here $\phi_{ji} = \phi_{ij}^{-1}$.
\end{list} 
We need to prove that both the obstruction for gluing the
$\mathcal{O}^{\prime}_{i}$  {\em and} the obstruction for lifting
$\mathcal{L}$ lie in $H^{2}(\mathcal{E}_{L})$. 

We have $\phi_{ji}\phi_{kj}\phi_{ik} =
\id_{\mathcal{O}_{\mathcal{X}}} + tD_{ijk}$ where $D_{ijk}$ is a \v{C}ech
2-cocycle of $\Theta_{X}$. This cycle represents the obstruction for gluing
the $\mathcal{O}^{\prime}_{i}$. We may assume $\mathcal{L}$ is given by
$F_{ij} \in \Gamma(U_{ij},
\mathcal{O}^{*}_{\mathcal{X}})$ satisfying the cocycle condition
$F_{ij}F_{jk}=F_{ik}$. Choose $F^{\prime}_{ij}
\in \Gamma(U_{ij}, (\mathcal{O}^{\prime}_{i})^{*})$ with
$\phi_{ij}(F^{\prime}_{ij}) = F^{\prime}_{ji}$ lifting the $F_{ij}$. Thus
$$
F^{\prime}_{ij} \phi_{ji}(F^{\prime}_{jk}) (F^{\prime}_{ik})^{-1} = 1 +
tg_{ijk}$$
for some $g_{ijk} \in \Gamma(U_{ijk}, \mathcal{O}_{X})$.

Since $e_{L}$ is locally split we may write $\mathcal{E}_{L}$ locally on
$U_{i}$ as $\mathcal{O}_{U_{i}} \oplus  \Theta_{U_{i}}$. The gluing is
determined (dually) by the extension class in $H^{1}(\Omega^{1}_{X})$;
$(g_{i}, D_{i}) \in \Gamma(U_{i},
\mathcal{E}_{L})$ and $ (g_{j}, D_{j}) \in
\Gamma(U_{j}, \mathcal{E}_{L})$ are equal on $U_{ij}$ iff $D_{i} = D_{j}$
and $g_{j} - g_{i} = D_{i}(f_{ij})/f_{ij}$. Now copy the proof of
\cite[Theorem 3.3.11 (ii)]{ser:def} to show that $(g_{ijk},D_{ijk})$
represents the obstruction in $\mathcal{E}_{L}$.

(iv)  This follows from considering commutative diagrams like \ref{B}. 

(v) This follows from a theorem of Grothendieck, 
\cite[Theorem 2.5.13]{ser:def}, and the proof of
\cite[Theorem 2.5.13]{ser:def}). 
\end{proof}

\section{$T^1_{A_{\mathcal{K}}}$ and $T^2_{A_{\mathcal{K}}}$ for manifolds}
\label{tiA}

We recall the description in \cite{ac:cot} of the multi-graded pieces of
$T^i_{A_\mathcal{K}}$ for any complex $\mathcal{K}$.  We will often denote
$T^i_\mathbf{c}(\mathcal{K}): = T^1_{A_\mathcal{K},\mathbf{c}}$ for
$\mathbf{c} \in \mathbb{Z}^{n+1}$. If $b \subseteq [n]$ let
\vspace{-0.5ex}
\[
U_{b} = U_{b}(\mathcal{K}) :=\{f \in \mathcal{K} : f\cup b \not\in
\mathcal{K}\}
\vspace{-0.5ex}
\] 
and
\[
\widetilde{U}_{b} = \widetilde{U}_{b}(\mathcal{K}) := \{f \in \mathcal{K} :
(f\cup b)\setminus \{v\} \not\in \mathcal{K} \text{ for some } v\in b\}
\subseteq U_{b}\,.
\]
Notice that $U_b=\widetilde{U}_{b}=\mathcal{K}$ unless $
\partial b$ is a subcomplex of $\mathcal{K}$. Moreover, if $
\partial b\subseteq \mathcal{K}$, then with
$L_{b}:=\bigcap_{b^\prime\subset b}\link(b^\prime,\mathcal{K})$ we have
\begin{equation*}
\mathcal{K}\setminus U_b = 
\begin{cases}  \emptyset \\
\overline{\st}(b)
\end{cases}
\hspace{0.0em}\mbox{and}\hspace{0.7em}
\mathcal{K}\setminus \widetilde{U}_{b} = 
\begin{cases}
\partial b \ast L_b& \text{if $b$ is a non-face},\\
(
\partial b \ast L_b) \cup \overline{\st}(b)& \text{if $b$ is a face}.
\end{cases}
\vspace{1ex}
\end{equation*}

\begin{theorem}\emph{(\cite[Theorem 13]{ac:cot})}\label{topopen} The
homogeneous pieces in degree $\mathbf{c}=\mathbf{a}-\mathbf{b}$ (with
disjoint supports $a$ and $b$) of the cotangent cohomology of the
Stanley-Reisner ring $A_\mathcal{K}$ vanish unless $a\in \mathcal{K}$,
$\mathbf{b} \in \{0,1\}^{n+1}$, 
$b\subseteq [\link(a)]$ and $b\ne
\emptyset$.  If these conditions are satisfied, we have isomorphisms
\[
T^i_\mathbf{c}(\mathcal{K})
\;\simeq\; 
H^{i-1}\big(\langle U_{b}(\link(a,\mathcal{K}))\rangle, \,
\langle \widetilde{U}_{b}(\link(a,\mathcal{K}))\rangle,\,\mathbb{C}\big)
\;\text{ for } i=1,2
\]
unless $b$ consists of a single vertex.  If $b$ consists of only one
vertex, then the above formulae become true if we use the reduced
cohomology instead.
\end{theorem} 

Since $T^i_\mathbf{c}(\mathcal{K})$ depends only on the supports $a$ and
$b$ we will often denote it $T^i_{a-b}(\mathcal{K})$. We will now apply the
result to combinatorial manifolds. We may reduce the computation to the
$a=\emptyset$ case by
\begin{proposition}\emph{(\cite[Proposition 11]{ac:cot})}\label{aempty} If
$\,b\subseteq [\link(a)]$, then the map $f\mapsto f\setminus a$  
induces isomorphisms $T^i_{\emptyset -b}(\link(a, \mathcal{K})) \simeq
T^i_{a-b}(\mathcal{K})$ for $i=1,2$.
\end{proposition}

\begin{lemma}\label{con} If $\mathcal{K}$ is a manifold and $b \neq
\emptyset$, then $U_{b}(\mathcal{K})$ is never empty and $\langle
U_{b}(\mathcal{K})\rangle$ is connected. Thus $$
\dim_{k}T^1_{\emptyset-b}(\mathcal{K}) =
\begin{cases}
        1 & \text{if $\widetilde{U}_{b}(\mathcal{K}) =
\emptyset$ and $|b| \ge 2$},
\\
        0& \text{otherwise}\, .
\end{cases}$$
\end{lemma}
\begin{proof} Set $U:=U_{b}(\mathcal{K})$.  If $b\not\in
\mathcal{K}$, then $\emptyset\in U$. Thus $U$ is non-empty and $\langle
U\rangle$ is a cone, so connected. If $b \in \mathcal{K}$ and
$U=\emptyset$, then $\mathcal{K} = \overline{\st}(b)$; i.e. a ball.  This
contradicts $\mathcal{K}$ being without boundary. If $b \in \mathcal{K}$
then $|\mathcal{K}|
\setminus \langle U\rangle = |\overline{\st}(b)|$, in particular
contractible.  Since $\mathcal{K}$ is a manifold, $\langle U\rangle$ is
connected.
\end{proof} 

\begin{remark} One can use the results of \cite{ac:cot} to compute the
$T^i$ also when $\mathcal{K}$ has boundary.  In this case though the
$U_{b}$ may not be connected if $b$ is a face and we do not get as nice
formulae as we do in the non-boundary case.
\end{remark}

\begin{definition} \label{bk} Define $\mathcal{B}(\mathcal{K})$ to be the
set of $b\subseteq [\mathcal{K}]$, $|b| \ge 2$, with the properties 
\begin{list}{\textup{(\roman{temp})}}{\usecounter{temp}}
\item $\mathcal{K}=L\ast 
\partial b$ where $|L|$ is a $(n-|b|+1)$-sphere if $b\not\in \mathcal{K}$,
\item $\mathcal{K}=L\ast 
\partial b \cup 
\partial L \ast
\bar{b}$ where $|L|$ is a $(n-|b|+1)$-ball if $b\in \mathcal{K}$.
\end{list}
\end{definition}
 Note that if $\mathcal{K}$ is not a sphere, then $\mathcal{B}(\mathcal{K})
= \emptyset$.
\begin{lemma}\label{U_tilde_empty} If $\mathcal{K}$ is an $n$-manifold and
$|b|\ge 2$ then $\widetilde{U}_{b}(\mathcal{K}) = \emptyset$ iff $b \in
\mathcal{B}(\mathcal{K})$.
\end{lemma}

\begin{proof} If $b \notin \mathcal{K}$ then
$\widetilde{U}_{b}(\mathcal{K}) = \emptyset$ means that 
$\mathcal{K} = L_b  \ast 
\partial b$. If $F$ is a facet of $
\partial b$, then $L_b = \link(F,\mathcal{K})$ is a sphere. If $b \in
\mathcal{K}$ then $\widetilde{U}_{b}(\mathcal{K}) = \emptyset$ means that
$\mathcal{K} = (L_b \ast \partial b) \cup \overline{\st}(b)$, i.e.
$\mathcal{K}\setminus \st(b) = L_b \ast 
\partial b$.  Now $\mathcal{K}\setminus \st(b)$ is a manifold with boundary
and $\partial b$ is in this boundary.  If $F$ is a facet of $
\partial b$, then $L_b = \link(F,\mathcal{K}\setminus \st(b))$ and
therefore a ball.
\end{proof}

We may add up these results to get a description of the whole
$T^1_{A_\mathcal{K}}$. 
\begin{theorem}\label{whole_t1} If $\mathcal{K}$ is a manifold and
$\mathbf{c}=\mathbf{a}-\mathbf{b}$ (with disjoint supports $a$ and $b$)
then 
$$
\dim_{k}T^1_{A_\mathcal{K},\mathbf{c}} =
\begin{cases}
        1 & \text{if $a\in \mathcal{K}$ and $b \in  \mathcal{B}(\link(a,
\mathcal{K}))$},
\\
        0& \text{otherwise}\, .
\end{cases}$$
A basis for $T^1_{A_\mathcal{K}}$ may be explicitly described: if $\phi \in
T^1_{A_\mathcal{K},\mathbf{c}} \neq 0$ and $x_{p} \in I_{\mathcal{K}}$ then
$\phi(x_{p}) = x^{\mathbf{a}}x_{p\setminus b}$ if $b \subseteq p$ and $0$
otherwise.
\end{theorem}
\begin{proof} This follows from Lemma~\ref{con}, Proposition~\ref{aempty}
and Lemma~\ref{U_tilde_empty}.
\end{proof}
\begin{remark} The case where $b$ is not a face corresponds to the notion
of {\it stellar exchange} defined in \cite{pac:plh}. (See also
\cite{vir:lec}.)  Assume $\mathcal{K}$ is a complex with a non-empty face
$a$ such that $\link(a,\mathcal{K})=
\partial b \ast L$ for some non-empty set $b$ and  $b$ is not a face of
$\link(a,\mathcal{K})$. We can now make a new complex
$\flip_{a,b}(\mathcal{K})$ by removing $\overline{\st}(a)= 
\partial b \ast \bar{a} \ast L$ and replacing it with $
\partial a \ast \bar{b} \ast L$, $$
\flip_{a,b}(\mathcal{K}):= (\mathcal{K}\setminus (
\partial b \ast \bar{a}
\ast L)) \cup 
\partial a \ast \bar{b} \ast L \, .$$
If $|b| = 1$, that is if $b$ is a new vertex, then
$\flip_{a,b}(\mathcal{K})$ is just the ordinary result of starring $b$ at
$a$. We see from Theorem~\ref{whole_t1} that if $a$ is not empty and $b$ is
not a face, then $a-b$ contributes to $T^1$ exactly when we can construct
$\flip_{a,b}(\mathcal{K})$.
\end{remark}

In dimensions $0$, $1$ and $2$ we may classify all the manifolds with
$\mathcal{B}(\mathcal{K}) \neq \emptyset$.  We use the notation of
Section~\ref{simp}. 
If $X$ is finite set, let $\mathcal{P}_{n}(X) \subseteq 2^{X}$ be the set
of subsets $Y$ with $\lvert Y \rvert = n$. Set $\mathcal{P}_{\ge n}(X) =
\bigcup_{r \ge n} \mathcal{P}_{r}(X)$.

\begin{table}
\centering
\begin{tabular}{|l|l|c|}    \hline $\mathcal{K}$ &
$\mathcal{B}(\mathcal{K})$ & $\lvert \mathcal{B}(\mathcal{K}) \rvert$ \\
\hline \hline  $
\partial \bigtriangleup_1$ & $\{[\mathcal{K}]\}$ &  $1$ \\
\hline \hline $
\partial \bigtriangleup_2$ &  $\mathcal{P}_{\ge 2}([\mathcal{K}]) $ &  $4$
\\
\hline $E_{4} = \mathcal{K}_{1} \ast \mathcal{K}_{2}$,  $\mathcal{K}_{i}
=\partial \bigtriangleup_1
 $& 
$\{[\mathcal{K}_{1}], [\mathcal{K}_{2}]\}$ & 
$2$
\\
\hline \hline $
\partial \bigtriangleup_3$   & $\mathcal{P}_{\ge 2}([\mathcal{K}]) $ & 
$11$ \\
\hline $\Sigma E_3 = 
\partial \bigtriangleup_1 \ast 
\partial\bigtriangleup_2$ & $\mathcal{B}(
\partial\bigtriangleup_1) \cup
\mathcal{B}(
\partial\bigtriangleup_2)$ & $5$
\\
\hline  $\Sigma E_4 = 
\partial\bigtriangleup_1 \ast 
E_{4}$  &  $\mathcal{B}(
\partial\bigtriangleup_1) \cup
\mathcal{B}( E_{4})$ & $3$\\
\hline  $\Sigma E_n = 
\partial\bigtriangleup_1 \ast  E_{n}, n\ge 5$ & $\{[
\partial\bigtriangleup_1]\}$ & $1$ \\
\hline  $
\partial C(n,3) , n\ge 6$ & 
$\{[
\partial\bigtriangleup_1]\}$ & $1$ \\
\hline
\end{tabular}
\caption{Manifolds $\mathcal{K}$ with $\dim \mathcal{K} \le 2$ and
$\mathcal{B}(\mathcal{K}) \neq \emptyset$.}
\label{lowdim}
\end{table}

\begin{proposition}\label{t1_list} If $\mathcal{K}$ is a manifold and $\dim
\mathcal{K}\le 2$, then 
$\mathcal{B}(\mathcal{K}) \neq \emptyset$ if and only if $\mathcal{K}$ is
one of the triangulations in Table~\ref{lowdim}.
\end{proposition}

We are not able to get so precise results for $T^2$, but for oriented
manifolds and especially spheres, $T^2$ is reasonably computable. Again it
is enough to compute the case $a=\emptyset$ and then use these results on
$\link(a)$ in the general case.

\begin{proposition}\label{t2_empty} If $\mathcal{K}$ is an $n$-manifold
then $T^2_{\emptyset-b}=0$ unless $\partial b \subset \mathcal{K}$.  If $
\partial b \subset \mathcal{K}$ and $L_{b} =
\cap_{b^\prime\subset b}\link(b^\prime,\mathcal{K})$, then
$T^2_{\emptyset-b}$ may be computed as follows:
\begin{list}{\textup{(\roman{temp})}}{\usecounter{temp}}
\item If $b\not\in \mathcal{K}$, then $T^2_{\emptyset-b}\simeq
\widetilde{H}^{0}(|\mathcal{K}|\setminus |
\partial b \ast L_{b}|,k)$.  If 
$|\mathcal{K}|$ is a sphere, then $T^2_{\emptyset-b}\simeq
\widetilde{H}_{n-|b|}(L_{b},k)$.
\item If $b\in \mathcal{K}$, then $T^2_{\emptyset-b}\simeq
H^{1}(|\mathcal{K}|\setminus |\overline{\st}(b)|,|\mathcal{K}|\setminus |(
\partial b \ast L_b) \cup
\overline{\st}(b)|,k)$. If $b$ is a vertex and $\mathcal{K}$ is oriented,
then $T^2_{\emptyset-b}\simeq \widetilde{H}_{n-1}(\mathcal{K},k)$. If $|b|
\ge 2$ and $\mathcal{K}$ is oriented, then $T^2_{\emptyset-b}=0$ if
$T^1_{\emptyset-b}\ne 0$.  If $T^1_{\emptyset-b}= 0$ then there is an exact
sequence $$
0 \rightarrow
\widetilde{H}_{n-|b|}(\link(b),k) \rightarrow
\widetilde{H}_{n-|b|}(L_{b},k) \rightarrow T^2_{\emptyset-b}
\rightarrow 0\, .$$
In particular $\dim T^2_{\emptyset-b} = \max \{\dim
\widetilde{H}_{n-|b|}(L_{b},k) - 1, 0\}$.
\end{list} These results are true even when the degree $n-|b|=-1$ with the
convention $\widetilde{H}_{-1}(\emptyset) = k$.  If $b^\prime$ is a facet
of $
\partial b$, then $\widetilde{H}_{n-|b|}(L_{b})$ may be computed as
$\widetilde{H}^{0}(\link(b^\prime)\setminus L_{b})$.
\end{proposition}

\begin{proof} By Theorem~\ref{topopen} we have $T^2_{\emptyset-b}$
isomorphic with $H^{1}(\langle U_{b}\rangle, \langle
\widetilde{U}_{b}\rangle)$.  If $b\not\in \mathcal{K}$, then $\emptyset \in
U_{b}$, so $\langle U_{b}\rangle$ is a cone.  Thus $H^{1}(\langle
U_{b}\rangle, \langle
\widetilde{U}_{b}\rangle)\simeq
\widetilde{H}^{0}(|\mathcal{K}|\setminus |
\partial b \ast L_{b}|,k)$. If $\mathcal{K}$ is a sphere, then by Alexander
duality $\widetilde{H}^{0}(|\mathcal{K}|\setminus |
\partial b \ast L_{b}|) \simeq
\widetilde{H}_{n-1}(
\partial b \ast L_b)$.  Now $|
\partial b|$ is homeomorphic to $S^{|b|-2}$, so $|
\partial b \ast L|$ is homeomorphic to the $(|b|-1)$-fold suspension of
$|L|$.  Thus $\widetilde{H}_{n-1}(
\partial b \ast L_b) \simeq
\widetilde{H}_{n-|b|}(L_{b})$. 

If $|b|=1$, then $\widetilde{U}_{b}=\emptyset$.  If $\mathcal{K}$ is
oriented then by duality $T^2_{\emptyset-b} \simeq H_{n-1}(\mathcal{K},
\overline{\st}(b)) \simeq \widetilde{H}_{n-1}(\mathcal{K})$.

If $b \in \mathcal{K}$ and $|b|\ge 2$ use first duality to get
$T^2_{\emptyset-b}\simeq H_{n-1}(\partial b \ast L_b \cup
\overline{\st}(b),
\overline{\st}(b))$.  Since $|b|\ge 2$, if we excise $\st(b)$, we achieve
an isomorphism with $H_{n-1}(\partial b \ast L_b, 
\partial b
\ast \link(b))$.  Again, because $|b|\ge 2$, $T^1_{\emptyset-b} \simeq
H^{0}(\langle U_{b}\rangle, \langle \widetilde{U}_{b}\rangle) \simeq H_{n}(
\partial b \ast L_b, 
\partial b \ast \link(b))$. Now $
\partial b \ast \link(b)$ is an $(n-1)$-sphere, so if $T^1_{\emptyset-b}=0$
we get an exact sequence $$
0 \rightarrow H_{n-1}(
\partial b \ast \link(b)) \rightarrow H_{n-1}(
\partial b \ast L_b ) \rightarrow T^2_{\emptyset-b}
\rightarrow 0\, .$$
The suspension argument gives the exact sequence in the statement.

If $T^1_{\emptyset-b}\ne 0$, then $\mathcal{K}=
\partial b \ast L_b \cup
\overline{\st}(b)$ by Lemma~\ref{U_tilde_empty}.  In particular $L_b =
\link(b^\prime) \approx S^{n-|b|+1}$ for all maximal $b^\prime \subset b$
and $(\partial b \ast L_b \cup \overline{\st}(b),
\overline{\st}(b)) \approx (S^n, B^n)$.

The last statement follows from Alexander duality on the $(n-|b|+1)$-sphere
$\link(b^\prime)$.
\end{proof}

\begin{remark} For 2-dimensional spheres an analysis yields the list of
unobstructed rings in \cite[Corollary 2.5]{iso:tor}.
\end{remark}

\section{$T^1_{{\mathbb P}(\mathcal{K})}$ and $T^2_{{\mathbb
P}(\mathcal{K})}$ for manifolds}\label{c_def}

We recall from \cite{ac:cot} the description of the derivations of
$A_\mathcal{K}$.
\begin{proposition}\emph{(\cite[Corollary 10]{ac:cot})}\label{T0}
$\,T^0_{A_\mathcal{K}} = \bigoplus_{v=0}^n
\mathfrak{a}_v\, 
\partial / 
\partial x_v$ where $\mathfrak{a}_v$ is the ideal of $A_\mathcal{K}$
generated by the monomials $x_{a}$ with $\overline{\st}(a,\mathcal{K})
\subseteq
\overline{\st}(v,\mathcal{K})$. In particular, $T^0_{A_\mathcal{K}}$ is
generated, as a module, by $x_v \, 
\partial / 
\partial x_v$ if and only if every non-maximal $a\in \mathcal{K}$ is
properly contained in at least two different faces.
\end{proposition} 
Certainly the criteria of the second statement is met by manifolds (without
boundary).  We may exploit this to construct an ``Euler sequence" for
${\mathbb P}(\mathcal{K})$. Let $y_j^{(i)}=x_j/x_i$ be coordinates for
$D_{+}(x_i)$ and set $\delta_j^{(i)} = y_j^{(i)} \, \partial/ \partial
y_j^{(i)}$.  By the global sections $\delta_i = x_i \, 
\partial/ 
\partial x_i$ we mean the \v Cech global sections $$
\delta_i =(\delta_i^{(0)},\dots ,\delta_i^{(i-1)}, -\sum_{j\ne i}
\delta_j^{(i)}, \delta_i^{(i+1)}, \dots, \delta_i^{(n)}) $$
which are subject to the relation $\sum_{i=0}^n \delta_{i} = 0$.

Let $S_i = \mathbb{P}(\overline{\st}(\{i\},\mathcal{K})) \subset
\mathbb{P}(\mathcal{K})$ where we view $S_i$ as embedded in
$\mathbb{P}^{n}$, i.e.\ $I_{S_i}$ contains all $x_j$ with $\{j\} \cup
\{i\} \not\in \mathcal{K}$.

\begin{theorem}\label{euler} If $\mathcal{K}$ is a manifold, then there is
an exact sequence of sheaves $$
0 \to \mathcal{O}_{\mathbb{P}(\mathcal{K})} \to
\bigoplus_{i=0}^n
\mathcal{O}_{S_i} \to \Theta_{\mathbb{P}(\mathcal{K})} \to 0\, .$$
The cohomology of $\Theta_{\mathbb{P}(\mathcal{K})}$ is given by
$H^p(\mathbb{P}(\mathcal{K}), \Theta_{\mathbb{P}(\mathcal{K})}) \simeq
H^{p+1}(\mathcal{K},\mathbb{C})$ if $p\ge 1$ and the exact sequence $$
0 \to \mathbb{C}^{n} \to H^0(\mathbb{P}(\mathcal{K}),
\Theta_{\mathbb{P}(\mathcal{K})}) \to H^1(\mathcal{K},\mathbb{C}) \to 0\,
.$$
\end{theorem}

\begin{proof} By Proposition~\ref{T0}, $\Theta_{\mathbb{P}(\mathcal{K})}$
is generated by the global sections $\delta_{i}$.  This gives a surjection
$\mathcal{O}_{\mathbb{P}(\mathcal{K})}^n \to
\Theta_{\mathbb{P}(\mathcal{K})}$. The annihilator of $\delta_i$ is the
ideal sheaf associated to $\Ann x_i \subseteq A_\mathcal{K}$. Clearly $\Ann
x_i + I_\mathcal{K}$ is the Stanley-Reisner ideal of
$\overline{\st}(\{i\},\mathcal{K})$.

The natural homomorphisms $A_\mathcal{K} \to A_\mathcal{K}/\Ann x_i$ add up
to an injection $A_\mathcal{K} \to \bigoplus A_\mathcal{K}/\Ann x_i$ since
every non-empty $f\in
\mathcal{K}$ is in some $\overline{\st}(\{i\})$. This gives the exact
sequence. Applying cohomology to this sequence yields the second statement.
Indeed, $\overline{\st}(\{i\})$ is contractible so the isomorphisms follow
from Theorem~\ref{O_cohom}.
\end{proof}

Let $B_i = \mathbb{P}(\mathcal{K} \setminus \st(\{i\},\mathcal{K})) \subset
\mathbb{P}(\mathcal{K})$ where we view $B_i$ as embedded in
$\mathbb{P}^{n}$, i.e.\ $I_{B_i} = I_\mathcal{K} + \langle x_i \rangle$.

\begin{proposition}\label{ob} If $\mathcal{K}$ is a manifold, then in the
exact sequence $$
0 \to \Theta_{\mathbb{P}(\mathcal{K})} \stackrel{\gamma}{\to}
\Theta_{\mathbb{P}^{n}}\otimes \mathcal{O}_{\mathbb{P}(\mathcal{K})} \to
{\mathcal N}_{{\mathbb P}(\mathcal{K})} \stackrel{\delta}{\to} {\mathcal
T}_{\mathbb{P}(\mathcal{K})}^1 \to 0$$
we have $\Ker(\delta) = \Coker(\gamma) \simeq
\bigoplus_{i=0}^n \mathcal{O}_{B_i}(1)$.
\end{proposition}

\begin{proof} By Theorem~\ref{euler} there is a commutative diagram of
Euler sequences with exact rows $$
\begin{CD} 0 @>>> \mathcal{O}_{\mathbb{P}(\mathcal{K})} @>>>
\bigoplus_{i=0}^n
\mathcal{O}_{S_i} @>>> \Theta_{\mathbb{P}(\mathcal{K})} @>>> 0\\
@. @VV\alpha V @VV\beta V @VV\gamma V @.\\
0 @>>> \mathcal{O}_{\mathbb{P}(\mathcal{K})} @>>> \bigoplus_{i=0}^n
\mathcal{O}_{\mathbb{P}(\mathcal{K})}(1) @>>>
\Theta_{\mathbb{P}^{n}}\otimes
\mathcal{O}_{\mathbb{P}(\mathcal{K})} @>>> 0
\end{CD}$$
where $\alpha$ is the identity and $\beta$ is induced from multiplication
with the $x_i$.  Thus the cokernel of $\gamma$ equals the cokernel of
$\beta$ which is clearly $\bigoplus_{i=0}^n
\mathcal{O}_{B_i}(1)$.
\end{proof}

For the local Hilbert functor $\Def_{{\mathbb P}(\mathcal{K})/{\mathbb
P}^{n}}$ we have the following result which we will also need in the
sequel.
\begin{proposition}\label{hilb} If $\mathcal{K}$ is a simplicial complex
then
\begin{list}{\textup{(\roman{temp})}}{\usecounter{temp}}
\item $H^0({\mathbb P}(\mathcal{K}),{\mathcal N}_{{\mathbb
P}(\mathcal{K})/{\mathbb P}^{n}}) \simeq
\Hom_P(I_\mathcal{K},A_\mathcal{K})_0$,
\item $T^2_{{\mathbb P}(\mathcal{K})/{\mathbb P}^{n}} \simeq
T^2_{A_\mathcal{K},0}$ and $T^2_{A_\mathcal{K},0} \to
H^0(\mathbb{P}(\mathcal{K}), {\mathcal T}_{\mathbb{P}(\mathcal{K})}^2)$ is
injective.
\end{list}
\end{proposition}
\begin{proof} The first statement follows from Schlessinger's comparison
theorem, see \cite{ps:hil} or \cite[Theorem 9.1]{ser:top}. For the second
statement, a close look at Kleppe's proof of the comparison theorem (see
\cite[3]{kle:def}) shows that if $H^0_{\mathfrak{m}}(A)=0$ and both
$H^1_{\mathfrak{m}}(A)$ and $H^2_{\mathfrak{m}}(A)$ vanish in positive
degrees, then $(T^2_{A})_0
\simeq T^2_{\Proj A/{\mathbb P}^r}$. Now apply Theorem~\ref{loco}. The
injectivity statement is \cite[Theorem 15]{ac:cot}.
\end{proof}

We are now able to describe the $T^i_{{\mathbb P}(\mathcal{K})}$.

\begin{theorem}\label{T_man} If $\mathcal{K}$ is a manifold then
\begin{list}{\textup{(\roman{temp})}}{\usecounter{temp}}
\item $H^0(\mathbb{P}(\mathcal{K}), {\mathcal
T}_{\mathbb{P}(\mathcal{K})}^1) \simeq T^1_{A_\mathcal{K},0}$.
\item $H^1(\mathbb{P}(\mathcal{K}), {\mathcal
T}_{\mathbb{P}(\mathcal{K})}^1) = 0$.
\item There are exact sequences
\begin{align*} 0 &\to
H^1(\mathbb{P}(\mathcal{K}),\Theta_{\mathbb{P}(\mathcal{K})}) \to
T^1_{{\mathbb P}(\mathcal{K})}\to H^0(\mathbb{P}(\mathcal{K}), {\mathcal
T}_{\mathbb{P}(\mathcal{K})}^1) \to 0\\
0 &\to H^2(\mathbb{P}(\mathcal{K}),\Theta_{\mathbb{P}(\mathcal{K})}) \to
T^2_{{\mathbb P}(\mathcal{K})} \to H^0(\mathbb{P}(\mathcal{K}), {\mathcal
T}_{\mathbb{P}(\mathcal{K})}^2)\, .
\end{align*}
\end{list}
\end{theorem}

\begin{proof} We have $H^i(\mathcal{O}_{B_i}(1)) = 0$ when $i\ge 1$ by
Theorem~\ref{O_cohom}, so the map $H^0({\mathcal N}_{{\mathbb
P}(\mathcal{K})})
\to H^0({\mathcal T}_{\mathbb{P}(\mathcal{K})}^1)$ is surjective and
$H^i({\mathcal N}_{{\mathbb P}(\mathcal{K})}) \simeq H^i({\mathcal
T}_{\mathbb{P}(\mathcal{K})}^1)$ when $i\ge 1$.  Since $H^0({\mathcal
N}_{{\mathbb P}(\mathcal{K})}) \simeq
\Hom_P(I_\mathcal{K},A_\mathcal{K})_0$ by Proposition~\ref{hilb}, the exact
sequence in Proposition~\ref{ob} 
yields (i).

Since $H^1({\mathcal N}_{{\mathbb P}(\mathcal{K})}) \simeq H^1({\mathcal
T}_{\mathbb{P}(\mathcal{K})}^1)$ is the kernel of $T^2_{{\mathbb
P}(\mathcal{K})/{\mathbb P}^{n}} \rightarrow H^0({\mathcal T}^2_{{\mathbb
P}(\mathcal{K})})$, (ii) follows from Proposition~\ref{hilb}. 

The exact sequences come from the edge exact sequences of the global-local
spectral sequence for $T^i_{{\mathbb P}(\mathcal{K})}$, see e.g.\
\cite[\S 4]{pal:def}.  The surjectivity in the first sequence follows from
the exactness of 
$$
T^1_{{\mathbb P}(\mathcal{K})}\to H^0(\mathbb{P}(\mathcal{K}), {\mathcal
T}_{\mathbb{P}(\mathcal{K})}^1) 
\stackrel{d_2}{\to}
H^2(\mathbb{P}(\mathcal{K}),\Theta_{\mathbb{P}(\mathcal{K})}) \, .$$ 
By Proposition~\ref{ob} this $d_2$ factors through $H^1(\bigoplus_{i=0}^n
\mathcal{O}_{B_i}(1)) = 0$, so it is the zero map.  This, together with
$H^1(\mathbb{P}(\mathcal{K}), {\mathcal T}_{\mathbb{P}(\mathcal{K})}^1) =
0$ yields the second exact sequence as well.
\end{proof}

We may use the analysis in section \ref{tiA} to find formulae for $T^{1}$
and $T^{2}$ for 
low dimensional $\mathcal{K}$. Let $f_i$ be the number of $i$-dimensional
faces of $\mathcal{K}$ and let $f_i^{(k)}$ be number of $i$-dimensional
faces with valency $k$. 

\begin{theorem}\label{low_dim} If $\mathcal{K}$ is a $2$-dimensional
manifold then 
\begin{align*}
\dim T^1_{{\mathbb P}(\mathcal{K})} &= 4f_0^{(3)} + 2f_0^{(4)} + f_1 +
h^2(\mathcal{K})\\
&= f_0 + 9 \chi(\mathcal{K}) + h^2(\mathcal{K}) + \sum_{k\ge
6}2(k-5)f_0^{(k)}\\
h^2(\Theta_{\mathbb{P}(\mathcal{K})}) &= 0 \text{ and } \dim
T^2_{A_\mathcal{K},0} =
\sum_{k\ge 6} \frac{1}{2} k(k-5)f_0^{(k)} \, .
\end{align*} 
\end{theorem} If $\dim \mathcal{K} =3$ set
\begin{align*} d_{3} &= \#\{v \in \mathcal{K}: \link(v) = 
\partial \Delta_3\}\\
e_{3} &= \#\{v \in \mathcal{K}: \link(v) = \Sigma E_3 \}\\
e_{4} &= \#\{v \in \mathcal{K}: \link(v) = \Sigma E_4 \}\\
e_{\ge 5} &= \#\{v \in \mathcal{K}: \link(v) = \Sigma E_n \text{ for some }
n\ge 5\}\\
c_{\ge 6} &= \#\{v \in \mathcal{K}: \link(v) = 
\partial C(n,3) \text{ for some } n\ge 6\}\, .
\end{align*} 
\begin{theorem}\label{low_dim3} If $\mathcal{K}$ is a $3$-dimensional
manifold then 
 $$
\dim T^1_{{\mathbb P}(\mathcal{K})} = 11d_{3} + 5e_{3} + 3e_{4} + e_{\ge 5}
+ c_{\ge 6} +5f_1^{(3)} + 2f_1^{(4)} + h^2(\mathcal{K})\, .$$
\end{theorem}
\begin{proof}[Proof of Theorem \ref{low_dim} and Theorem \ref{low_dim3}] By
Theorem~\ref{T_man} and Theorem~\ref{euler} we need only to find the
contribution from $T^1_{A_\mathcal{K},0}$. The $T^1_{a-b}$ that contribute
in degree $0$ have $0 < |a|
\le |b|$.  By Theorem~\ref{whole_t1}, if $T^1_{a-b}\ne 0$, then $\dim
\mathcal{K} - \dim a + 1 \ge |b|$. 
We must therefore have $\dim a \le \frac{1}{2}\dim \mathcal{K}$.

Except for the case $\dim \mathcal{K} = 3$, $|a|=2$ and $|b|=3$, there is a
unique ${\mathbf a}$ making $|{\mathbf a}| = |b|$.  In the exceptional case
$\link(a)$ equals $
\partial \Delta_2$ and there are two choices for ${\mathbf a}$.  Thus
$f_1^{(3)}$ contributes with $5$.  The formulae for $\dim T^1_{A,0}$ can
now be computed from Proposition~\ref{t1_list}.

The second formula when $\dim \mathcal{K}= 2$ follows from 
$$
6\chi(\mathcal{K}) = \sum _{k \ge 3}(6 -k)f_0^{(k)} \, .$$
The $T^2$ formula follows from Proposition~\ref{t2_empty}.
\end{proof}

Since $f_1$ contributes to $T^1$ when $\mathcal{K}$ is a surface, ${\mathbb
P}(\mathcal{K})$ is never rigid in this case.  Things are different in
dimension $3$.

\begin{corollary}\label{rigid} If $\mathcal{K}$ is a $3$-dimensional
manifold, then ${\mathbb P}(\mathcal{K})$ is rigid if $H^{2}(\mathcal{K}) =
0$ and all edges $e$ have $\nu(e)\ge 5$.
\end{corollary}

\begin{example} If $\mathcal{K}$ is the boundary complex of the regular
solid with Schl{\"a}fli symbol $\{3,3,5\}$, then ${\mathbb P}(\mathcal{K})$
is rigid in ${\mathbb P}^{119}$.
\end{example}

We cannot give formulas for $T^{2}$ in the $3$-dimensional case, but
Proposition~\ref{t2_empty} is a useful tool for computations. We illustrate
this with a $3$-dimensional example.
\begin{example} Consider the boundary of the $4$-dimensional cyclic
polytope with $8$ vertices $\partial C(8,4)$ (see \cite[4.7]{gr:con}).
There are $20$ facets:
\begin{gather*}\{i,i+1,i+2,i+3\}, \{i,i+1,i+3,i+4\} \text{ for $i = 0,
\dots , 8$}\\
\{i,i+1,i+4,i+5\} \text{ for $i = 0,
\dots , 4$}
\end{gather*} where addition is modulo $8$. The links of the vertices are
all boundaries of the cyclic polytope $C(7,3)$. We draw the link of $\{0\}$
in Figure~\ref{C(7,3)}. 
We will compute $T^2_{A_{
\partial C(8,4)},0}$ using the statements and notation of
Proposition~\ref{t2_empty}.
\begin{figure}
\centering
\includegraphics[]{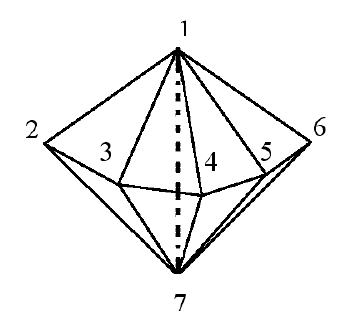}
\caption{The link of vertex $\{0\}$ which is $
\partial C(7,3)$. }
\label{C(7,3)}
\end{figure}

In dimension $3$, $T^{2}_{a-b} \ne 0$ with $|a| \le |b|$ implies that $\dim
a \le 1$. If $a$ is an edge then only the case  $\link(a) = E_{6}$
contributes to $T^2$ and the contribution may be computed as above (see
also \cite[Example 17]{ac:cot}). There are $8$ such edges, $\{i,i+1\}$, so
we get $8 \times 3 = 24$ basis elements this way.

If $a$ is a vertex we may assume by symmetry that $a = \{0\}$, so
$\link(a)$ is as drawn in Figure~\ref{C(7,3)}. We need to find the
different $b$ with the property  $T^{2}_{\emptyset - b}(\partial C(7,3))
\ne 0$.

Assume first $b$ is not a face. Thus  $T^{2}_{\emptyset - b}= 0$ if
$\partial b$ is not a sub-complex. If $\partial b$ is a sub-complex then
$T^{2}_{\emptyset - b} \simeq \widetilde{H}_{2-|b|}(L_{b},k)$. For $|b| =
2$, $L_{b}$ is empty or connected for all non-edges except $\{2,5\}$ and
$\{3,6\}$ for which $L_{b} = \{1\} \cup \{7\}$. For $|b| = 3$, $\partial b$
is a sub-complex 
for $\{1,3,7\}$, $\{1,5,7\}$ and $\{1,4,7\}$. Only $L_{\{1,4,7\}} =
\emptyset$.

Assume now $b$ is a face. If $b$ is a vertex then $T^{2}_{\emptyset - b}
\simeq H_{1}(
\partial C(7,3),k) = 0$. If $b = \{1,7\}$ then $T^{2}_{\emptyset - b} = 0$
since $T^{1}_{\emptyset - b} \ne 0$. For all other non-vertex faces we have
$\dim T^2_{\emptyset-b} = \max \{\dim
\widetilde{H}_{2-|b|}(L_{b},k) - 1, 0\}$. For this to be non-zero, $b$ must
be an edge and $L_{b}$ must have $3$ or more components. This happens only
for $\{1,4\}$ where $L_{b} = \{3\} \cup \{5\} \cup \{7\}$ and  $\{4,7\}$
where $L_{b} = \{1\} \cup \{3\} \cup \{5\}$.

Summing up we get a contribution to $T^{2}$ for $a = \{0\}$ when 
$b$ is  $\{2,5\}$, $\{3,6\}$,$\{1,4\}$,$\{4,7\}$ or $\{1,4,7\}$ and in each
case $\dim T^{2}_{a - b} = 1$. Thus all in all $\dim T^2_{A_{
\partial C(8,4)},0} = 24 + 8 \times 5 = 64$.
\end{example}

\section{Algebraic and non-algebraic deformations of ${\mathbb
P}(\mathcal{K})$}\label{algdef}

We consider now the functor $\Def_{{\mathbb P}(\mathcal{K})}^{a} :=
\Def_{({\mathbb P}(\mathcal{K}),L)}$, $L=\mathcal{O}_{{\mathbb
P}(\mathcal{K})}(1)$, of algebraic deformations. We will keep the notation
from Section~\ref{DXL} and \ref{c_def}. Recall that $S_i =
\mathbb{P}(\overline{\st}(\{i\},\mathcal{K}))$.

\begin{theorem} \label{Da} If $\mathcal{K}$ is a manifold then
$\mathcal{E}_{\mathcal{O}_{{\mathbb P}(\mathcal{K})}(1)} \simeq
\oplus_{i=0}^n \mathcal{O}_{S_i}$, in particular
$H^{i}(\mathcal{E}_{\mathcal{O}_{{\mathbb P}(\mathcal{K})}(1)} ) = 0$ for
$i \ge 1$. Thus $$
\Def_{{\mathbb P}(\mathcal{K})}^{a}(k[\epsilon]) \simeq
H^0(\mathbb{P}(\mathcal{K}), {\mathcal T}_{\mathbb{P}(\mathcal{K})}^1)
\simeq T^1_{A_\mathcal{K},0}$$ 
and $H^0(\mathbb{P}(\mathcal{K}), {\mathcal
T}_{\mathbb{P}(\mathcal{K})}^2)$ contains all obstructions for
$\Def_{{\mathbb P}(\mathcal{K})}^{a}$.
\end{theorem}
\begin{proof} We claim that the exact sequence in Theorem~\ref{euler}
represents the dual of $c(\mathcal{O}_{{\mathbb P}(\mathcal{K})}(1))$.
Indeed, from the proof of Theorem~\ref{DXL_theorem}, we see that
$\mathcal{E}_{L}$ is determined by being locally $\mathcal{O}_{U_{i}}
\oplus  \Theta_{U_{i}}$ with gluing ; $(g_{i}, D_{i}) \in \Gamma(U_{i},
\mathcal{E}_{L})$ and $ (g_{j}, D_{j}) \in
\Gamma(U_{j}, \mathcal{E}_{L})$ are equal on $U_{ij}$ iff $D_{i} = D_{j}$
and $g_{j} - g_{i} = D_{i}(f_{ij})/f_{ij}$. One checks that $\oplus_{i=0}^n
\mathcal{O}_{S_i}$ satisfies this when $f_{ij} = x_{j}/x_{i}$. The rest of
the statement follows from Theorem~\ref{DXL_theorem} and
Theorem~\ref{T_man}
\end{proof}

On the other hand we may consider the functor of locally trivial
deformations $\Def^{\prime}_{\mathbb{P}(\mathcal{K})}$. (See e.g.
\cite[1.1.2]{ser:def}.)
\begin{proposition} If $\mathcal{K}$ is a manifold then $$
\Def_{{\mathbb P}(\mathcal{K})}^{\prime}(k[\epsilon]) \simeq
H^1(\mathbb{P}(\mathcal{K}), {\Theta}_{\mathbb{P}(\mathcal{K})})
\simeq H^{2}(\mathcal{K},k)$$ 
and $H^2(\mathbb{P}(\mathcal{K}), {\Theta}_{\mathbb{P}(\mathcal{K})})
\simeq H^{3}(\mathcal{K},k)$ is an obstruction space for $\Def_{{\mathbb
P}(\mathcal{K})}^{\prime}$.
\end{proposition}
\begin{proof} This follows from Theorem~\ref{euler}.
\end{proof}

From now on let $\mathcal{K}$ be a $2$-manifold. If it is oriented then
$H^{2}(\mathcal{K},k)
\simeq k$ and $H^{3}(\mathcal{K},k) = 0$. Thus $\Def_{{\mathbb
P}(\mathcal{K})}^{\prime}$ has a smooth one dimensional versal base space.
If $k = \mathbb{C}$, since $H^{1}(\mathcal{E}_{\mathcal{O}_{{\mathbb
P}(\mathcal{K})}(1)} ) = 0$,  the fibers will consist of {\em
non-algebraic} deformations of the compact complex space $S = {\mathbb
P}_{\mathbb{C}}(\mathcal{K})$. We may describe them explicitly.

Let $y_{j}^{(i)} = x_{j}/x_{i}$ be local coordinates for $U_{i} =
D_{+}(x_{i})$.  As in Section~\ref{c_def} set $\delta_j^{(i)} = y_j^{(i)}
\,
\partial/ 
\partial y_j^{(i)}$. If $\{i,j\}$ is an edge set $U_{ij} = U_{i} \cap U_{j}
= D_{+}(x_{i}x_{j})$. If $\link(\{i,j\}) =
\{\{k\}, \{l\}\}$, then 
$$
U_{ij} = \Spec \mathbb{C}[y_{k}^{(i)},
y_{l}^{(i)},y_{j}^{(i)},(y_{j}^{(i)})^{-1}]/(y_{k}^{(i)} y_{l}^{(i)})$$
(see Section~\ref{sr_rings}) and the gluing is determined by $y_{j}^{(i)} =
x_{j}/x_{i}$.

We wish to understand the isomorphism $\mathbb{C} \simeq H^{2}(\mathcal{K},
\mathbb{C}) 
\simeq H^1({\Theta}_{S})$. If $\sigma$ is any oriented $2$-simplex of
$\mathcal{K}$, then the class of its dual $\sigma^{\ast}$ will be a
generator of $H^{2}(\mathcal{K},
\mathbb{C})
\simeq \mathbb{C}$. Assume $\sigma = \{i,j,k\}$ with $i < j < k$. One may
compute that the corresponding generator of $H^1({\Theta}_{S})$ is the
\v{C}ech cocycle $$
\delta_{\sigma} = \delta_k^{(i)}|_{U_{ij}} - \delta_j^{(i)}|_{U_{ik}} +
\delta_i^{(j)}|_{U_{jk}} \, .$$
The corresponding one parameter versal family over $\Delta = \{t \in
\mathbb{C} \, | \, |t| < 1\}$ is thus achieved by changing the gluing by 
\begin{align*} y_k^{(i)} &= (1-t)\frac{y_k^{(j)}}{y_i^{(j)}}   \quad 
\text{on $U_{ij}$}\\
y_j^{(i)} & = \frac{1}{(1-t)} \frac{y_j^{(k)}}{y_i^{(k)}}   \quad   
\text{on $U_{ik}$}\\
y_i^{(j)} & = (1-t)\frac{y_i^{(k)}}{y_j^{(k)}}   \quad   \text{on $U_{jk}$}
\end{align*} while all other identities remain the same. This defines a
family of complex spaces $\mathfrak{X} \to \Delta$.

We may describe this family in a way that generalizes the treatment of the
tetrahedron in \cite{fr:glo}. Let $P_{\sigma} \simeq \mathbb{P}^2$ be the
component of $S$ corresponding to $\sigma$, $S^{\prime} = 
\mathbb{P}(\mathcal{K} \setminus \sigma)$ and $D = S^{\prime} \cap 
P_{\sigma} \simeq \mathbb{P}(E_{3})$. Note that $S^{\prime}$ remains
unchanged by the new gluing since $\delta_{\sigma}| S^{\prime} = 0$. Of
course the restriction 
$\delta_{\sigma}|\mathbb{P}_{\sigma}$ is a coboundary and is the image of
$d=-\delta_j^{(i)}|_{U_{i}} + \delta_i^{(j)}|_{U_{j}}$.

By Theorem~\ref{euler} one sees that $H^{1}(E_{3})$ contributes to 
$H^{0}(D, \Theta_{D})$. This corresponds to a $\mathbb{C}^{\ast}$ action on
$D$ which is not induced by projective transformations of $\mathbb{P}^{2}$.
Now $d$ is a cocycle on $D$ and it's class in $H^{0}(\Theta_{D})$ generates
$H^{1}(E_{3})$. The corresponding family of automorphisms may be defined by
$\phi_{t}(x_{i}:x_{j}:0) = ((1-t)x_{i}:x_{j}:0)$ on the component $x_{k}=0$
and $\phi_{t} = 1$ on the other two components. We may regard $\phi_{t}$ as
an isomorphism 
$$
S^{\prime} \supset D \stackrel{\phi_{t}}{\simeq} D \subset P_{\sigma}\, .$$
We sum up the above in
\begin{proposition} If $\mathcal{K}$ is an oriented $2$-dimensional
manifold and $S=\mathbb{P}_{\mathbb{C}}(\mathcal{K})$ then the
$1$-dimensional versal locally trivial deformation $\mathfrak{X} \to
\Delta$ of $S$ has fibers $$
X_{t} \simeq S^{\prime} \sqcup P_{\sigma}/x \sim \phi_{t}(x) \, .$$
The fibers $X_{t}$, $t \neq 0$, are non-algebraic complex spaces.
\end{proposition}

We may compute $\Def^{a}_{S}$ when $\mathcal{K}$ is a $2$-dimensional
combinatorial manifold and all vertices $v$ have $\nu(v) \le 6$. Let $S =
\mathbb{P}(\mathcal{K})$. We start
by defining a set of coordinate functions corresponding dually to a basis
for $\Def^{a}_{S}(k[\epsilon])$. (See TheoremÊ~\ref{Da} and \cite[Example
18]{ac:cot}.) We need 
\begin{list}{}{\usecounter{temp}}
\item The variable $t_{i,j} = t_{j,i}$ for each edge $\{i,j\}$.
\item The $4$ variables $v_{i}, v_{i,j}, v_{i,k}, v_{i,l}$ for each vertex
$\{i\}$ with $\nu(\{i\}) = 3$ and $\{j\},\{k\}, \{l\}$ the vertices of
$\link(\{i\})$.
\item The $2$ variables $u_{i,i_{1}} = u_{i, i_{3}}$ and $u_{i,i_{2}} =
u_{i, i_{4}}$ for each vertex $\{i\}$ with $\nu(\{i\}) = 4$ and $\{i_{j},
i_{j+1}\}$ the edges of $\link(\{i\})$.
\end{list} Let $P_{S}$ be the polynomial $k$-algebra and
$\hat{P}_{S}$ the formal power series algebra in these variables.

For each vertex $\{ i_{0}\}$ with $\nu(\{i_{0}\}) = 6$, choose a cyclic
ordering of the vertices $\{i_{1}\},
\dotsc , \{i_{6}\}$ in the hexagon $\link(\{i_{0}\})$ so that $\{i_{j},
i_{j+1}\}$ are the edges of $\link(\{i_{0}\})$. Let $\mathfrak{G}_{i_{0}}$
be a set of $6$ power series in $\hat{P}_{S}$;
$\mathfrak{G}_{i_{0}} =
\{ g_{i_{0},i_{j}} | j = 1, \dotsc ,6 \}$. Set $$
\mathfrak{G} = \bigcup_{\nu(\{i\}) = 6} \mathfrak{G}_{i}$$
a set of $6f_0^{(6)}$ power series. Note that we do not assume $g_{i,j} =
g_{j,i}$ if both vertices have valency 6.

Let $\mathfrak{a}_{\mathfrak{G}_{i_{0}}} \subset \hat{P}_{S}$ be
the ideal generated by the $2 \times 2$ minors of $$
\begin{bmatrix} g_{i_{0},i_{1}} & g_{i_{0},i_{3}} &  g_{i_{0},i_{5}}\\
 g_{i_{0},i_{4}} &  g_{i_{0},i_{6}} & g_{i_{0},i_{2}}
\end{bmatrix}$$
and define the ideal 
\begin{equation}
\label{aS}
\mathfrak{a}_{\mathfrak{G}} = \sum_{\nu(\{i\}) = 6}
\mathfrak{a}_{\mathfrak{G}_{i_{0}}}\, .
\end{equation}
We set $\mathfrak{a}_{S} = \mathfrak{a}_{\mathfrak{G}}$ if all
$g_{ij} = t_{ij}$. Finally define the complete local $k$-algebra $
\hat{R}_{\mathfrak{G}} =
\hat{P}_{S}/\mathfrak{a}_{\mathfrak{G}}$. Denote the maximal
ideal of $\hat{R}_{\mathfrak{G}}$ by $\mathfrak{m}$.
\begin{theorem} \label{formal} If $\mathcal{K}$ is a $2$-dimensional
combinatorial manifold with $\nu(v) \le 6$ for all vertices, then we may
find $\mathfrak{G}$ as above with 
$$
g_{i,j} = t_{i,j} + \text{higher order terms}$$
such that $\Spec \hat{R}_{\mathfrak{G}}$ is a formal versal base space for
$\Def^{a}_{S}$. If $\nu(\{i\}) = 6$, $\{i, j\}$ is an edge and $\nu(\{j\})
\le 5$, then we may choose $g_{i,j} = t_{i,j}$.
\end{theorem}

\begin{example} If $\mathcal{K}$ is the suspension $\{\{0\}, \{7\}\} \ast
E_{6}$, then $\hat{P}_{S}$ is the power series ring in the $30$
variables $t_{0,j}$ for $j= 1,\dots ,6$ , $t_{7,j}$ for $j= 1,\dots ,6$,
$t_{i,i+1}$ for $i = 1,\dots ,6$, $u_{i,i+1} = u_{i, i-1}$ for  $i =
1,\dots ,6$ and $u_{i,0} = u_{i, 7}$ for $i = 1,\dots ,6$. The ideal
$\mathfrak{a}_{S}$ is generated by the $2 \times 2$ minors of
 $$
\begin{bmatrix} t_{0,1} & t_{0,3} &  t_{0,5}\\
  t_{0,4} &   t_{0,6} & t_{0,2}
\end{bmatrix} \text{  and   }
\begin{bmatrix} t_{7,1} & t_{7,3} &  t_{7,5}\\
  t_{7,4} &   t_{7,6} & t_{7,2}
\end{bmatrix} $$
and  $\hat{R}_{S}$ is the $26$ dimensional quotient ring.
\end{example}

We will prove the theorem using obstruction calculus. To do this we need to
know what the possible local deformations of each chart may look like. Let
$Z_{n} = \mathbb{A}(E_{n})$ and recall that $S$ is covered by $U_{i} \simeq
Z_{\nu(\{i\})}$. 

Index the vertices of $E_{n}$ cyclically by $1, 2, \dotsc, n$, all addition
is done modulo $n$, so that the edges of $E_{n}$ are $\{i , i+1\}$. The
Stanley-Reisner ideal of $Z_{n}$ for $n\ge 4$ is $I_{n} = (\{y_{i}y_{j}
: |j-i| \ge 2
\})$ in $k[y_{1},
\dotsc , y_{n}]$. 

The infinite dimensional $T^{1}_{Z_{n}}$ is computed in e.g. \cite{ac:cot}.
If $n\ge 5$ a basis may be represented by $\phi_{i}^{(k)}$, $k\ge 1$, 
which map $y_{i-1}y_{i+1} \mapsto y_{i}^{k}$ and all other generators of
the ideal to $0$. If $n=4$ then in addition we have $2$ basis elements,
both with two names, 
$\phi_{2}^{(0)} = \phi_{4}^{(0)}$ which maps $y_{1}y_{3} \mapsto 1,
y_{2}y_{4}
\mapsto 0$ and $\phi_{1}^{(0)} = \phi_{3}^{(0)}$ which maps $ y_{2}y_{4}
\mapsto 1, y_{1}y_{3} \mapsto 0$. Finally if $n=3$ we have the basis
$\phi_{i}^{(k)}: y_{1}y_{2}y_{3} \mapsto y_{i}^{k+1}$, $k \ge 0$ and
additionally $\phi^{(-1)}_{1} = \phi^{(-1)}_{2} = \phi^{(-1)}_{3}$ mapping
$y_{1}y_{2}y_{3} \mapsto 1$. We will denote the dual coordinate functions
in the symmetric algebra $\Sym(T^1_{Z_{n}})$ by $t_{i}^{(k)}$.

For $n =3, 4, 5, 6$ we will define a {\em normal form} for a 
deformation of $Z_{n}$. These will consist of 
a $k$-algebra $\mathcal{R}_{n}$ which is a quotient of the infinite
dimensional algebra of formal power series 
$k[[t_{i}^{(k)}]]$, by a finitely generated ideal $\mathfrak{a}_{n}$ and a
finite set of equations $\mathcal{I}_{n} \subset k[y_{1}, \dotsc ,
y_{n}][[t_{i}^{(k)}]]/\mathfrak{a}_{n}$.

\vspace{1ex}
\noindent 
$E_{3}$ (Hypersurface):  Define the algebra $\mathcal{R}_{3}: = k[[t_{i}^{(k)}]]$ for $i =
1, 2, 3$ and $k \ge -1$. Let $T_{i} =
\sum_{k=1}^{\infty} t_{i}^{(k)}y_{i}^{k}$ and $u=t^{(-1)}_{1} =
t^{(-1)}_{2} = t^{(-1)}_{3}$. The one equation $$
y_{1}y_{2}y_{3} + u + y_{1}(t_{1}^{(0)} + T_{1}) + y_{2}(t_{2}^{(0)} +
T_{2}) + y_{3}(t_{3}^{(0)} + T_{3})$$
is all that is in $\mathcal{I}_{3}$.

\vspace{1ex}
\noindent 
$E_{4}$ (Complete intersection): Define the algebra $\mathcal{R}_{4}: = k[[t_{i}^{(k)}]]$ for $i =
1, 2,3,4$ and $k \ge 0$. Let $T_{i} =
\sum_{k=1}^{\infty} t_{i}^{(k)}y_{i}^{k-1}$, $u=t^{(0)}_{2} = t^{(0)}_{4}$
and $v=t^{(0)}_{1} = t^{(0)}_{3}$. The two equations
\begin{gather*} y_{1}y_{3}+ u + y_{2}T_{2} + y_{4}T_{4}\\
y_{2}y_{4}+ v + y_{1}T_{1} + y_{3}T_{3}
\end{gather*} make up $\mathcal{I}_{4}$.

\vspace{1ex}
\noindent 
$E_{5}$ (Pfaffian): Define the algebra $\mathcal{R}_{5}: = k[[t_{i}^{(k)}]]$ for $i =
1,\dotsc ,5$ and $k \ge 1$. Let $T_{i} =
\sum_{k=1}^{\infty} t_{i}^{(k)}y_{i}^{k-1}$. The five equations $$
y_{i-1}y_{i+1} + y_{i}T_{i} - T_{i-2}T_{i+2}$$
for $i= 1, \dotsc , 5$ make up $\mathcal{I}_{5}$.

\vspace{1ex}
\noindent 
$E_{6}$ (First obstructed case):  Let $\mathfrak{a}_{6}$ be the ideal generated by the $2
\times 2$ minors of 
\begin{equation}
\label{obeq}
\begin{bmatrix} t_{1}^{(1)} & t_{3}^{(1)} & t_{5}^{(1)}\\
t_{4}^{(1)} & t_{6}^{(1)} & t_{2}^{(1)}
\end{bmatrix}\, .
\end{equation}
Define the algebra $\mathcal{R}_{6}: = k[[t_{i}^{(k)}]]/\mathfrak{a}_{6}$
for $i = 1,\dotsc ,6$ and $k \ge 1$. Let $s_{i} =
\sum_{k=2}^{\infty} t_{i}^{(k)}y_{i}^{k-2}$ and $S = \prod _{i=1}^6 s_{i}$.
Let $p(x)$ be a power series solution of the functional equation 
\begin{equation*}
\label{f-id} xp(x)^4 = p(x)+1
\end{equation*} and set $f = p(S)$ and $e = f/(f+2)$. The six equations
\begin{multline*}y_{i-1}y_{i+1}+(t^{(1)}_i+s_iy_i)y_i \\
+ s_{i+3}(e^2t^{(1)}_{i-2}t^{(1)}_{i+2}+efs_{i+2}t^{(1)}_{i-2}y_{i+2} +
eft^{(1)}_{i+2}s_{i-2}y_{i-2})\\
-s_{i-2}s_{i+2}(et^{(1)}_{i+3}+fs_{i+3}y_{i+3})^2 \\
+ e^2f^2s_{i-2}s_{i-1}s_{i+1}s_{i+2}s_{i+3}(t^{(1)}_{i})^2
\end{multline*} for $i = 1, \dotsc ,6$ and the three equations
\begin{multline*}y_iy_{i+3} + et^{(1)}_{i+1}t^{(1)}_{i+2} +
et^{(1)}_{i+2}s_{i+1}y_{i+1} + et^{(1)}_{i+1}s_{i+2}y_{i+2} +
fs_{i+1}s_{i+2}y_{i+1}y_{i+2}\\
 + et^{(1)}_{i-2}s_{i-1}y_{i-1} + et^{(1)}_{i-1}s_{i-2}y_{i-2} +
fs_{i-1}s_{i-2}y_{i-1}y_{i-2}\\
 - e^2f^2s_{i-2}s_{i-1}s_{i+1}s_{i+2}t^{(1)}_it^{(1)}_{i+3}
\end{multline*} for $i = 1,2,3$ make up $\mathcal{I}_{6}$. (See
\cite[4.3]{ste:deg} for a description of a similar family.)

\begin{proposition} \label{nf} For any $k$-algebra homomorphism
$\mathcal{R}_{n}
\to A$, for $n = 3,4,5,6$, where $A$ is an artinian local $k$-algebra and
almost all $t_{i}^{(k)} \mapsto 0$, the image of $\mathcal{I}_{n}$ in
$A[y_{1},
\dotsc ,y_{n}]$ defines a deformation $\mathcal{Z} \to \Spec(A)$ of
$Z_{n}$.
\end{proposition}

\begin{proof} We must prove that the relations among the generators of
$I_{Z_{n}}$ in $k[y_{1},
\dotsc , y_{n}]$ lift over $\mathcal{R}_{n}$ to relations among the
elements in $\mathcal{I}_{n}$. This is trivially true for $n=3,4$ and
easily checked for $n=5$. We will now prove it for $n=6$. 

To shorten notation set $t_{i} = t_{i}^{(1)}$. Let $F_{i,j}$ be the
equation in 
$\mathcal{I}_{6}$ lifting $y_{i}y_{j}$. The dihedral group $D_{6}$ acts on
everything by permuting indices. The action is generated by e.g. the cycle
$(1,2,3,4,5,6)$ and the reflection $(2,6)(3,5)$. In particular it acts on
$\mathcal{I}_{6}$

There are $16$ generators of the relation module for $I_{Z_{6}}$ and they
split into two $D_{6}$ orbits; the orbits of $y_{5}(y_{1}y_{3}) -
y_{1}(y_{3}y_{5})$ and $y_{6}(y_{1}y_{3}) - y_{1}(y_{3}y_{6})$. Using the
$D_{6}$ symmetry it is enough to give liftings of these 2 relations and one
checks that the following two expressions are such liftings:
\begin{multline*}(y_5+ef^3s_2s_3s_4s_1s_6t_5)F_{1,3} 
-(y_1+ef^3s_2s_3s_4s_5s_6t_1)F_{3,5} \\ 
+ s_4s_6(eft_5+f^2s_5y_5)F_{4,6} -s_2s_6(eft_1+f^2s_1y_1)F_{2,6}\\
 -(et_4+fs_4y_4)F_{1,4}+ (et_2+fs_2y_2)F_{2,5}
\end{multline*}
\begin{multline*} y_6F_{1,3}+ef^2s_2s_3s_4s_5t_4F_{2,4} 
-(eft_2+f^2s_2y_2)s_3s_4s_5F_{3,5}\\
 -efs_4s_5t_6F_{4,6} + et_4s_5F_{1,5} -(ef^{-1}t_2+s_2y_2)F_{2,6} \\
+ s_{4}(et_5+fs_5y_5)F_{1,4} 
-y_1F_{3,6} \, .
\end{multline*} 
These equations and relations were originally conjectured after using Maple
to lift equations and relations to degree 19.
\end{proof} 

\begin{definition} An infinitesimal  deformation $\mathcal{Z} \to
\Spec(A)$ of $Z_{n}$ is in {\em normal form} if it is induced in the above
sense by $(\mathcal{R}_{n}, \mathcal{I}_{n})$, i.e. there exists a
$k$-algebra homomorphism $\mathcal{R}_{n} \to A$ where almost all
$t_{i}^{(k)} \mapsto 0$ and $I_{\mathcal{Z}} \subset A[y_{1}, \dotsc
,y_{n}]$ is generated by the image of $\mathcal{I}_{n}$.
\end{definition}

\begin{proof}[Proof of Theorem~\ref{formal}] We will construct by induction
Cartesian diagrams of deformations of $S$
\begin{equation}
\label{cart}
\begin{CD} X_{n} @>>> X_{n+1} \\
@VVV @VVV \\
\Spec R_{n} @>>> \Spec R_{n+1}
\end{CD}
\end{equation} 
where the $R_{n}$ are local artinian quotients of $P_{S}$ with
$R_{n} \simeq R_{n+1}/\mathfrak{m}^{n+1}$, $\mathfrak{m}$ is the maximal
ideal of $\hat{P}_{S}$, and $\hat{R} = \lim R_{n}$ is as in the
theorem. Set first $R_{0} = k$ and $R_{1} =
P_{S}/\mathfrak{m}^{2}$. Thus the Kodaira-Spencer map will be
surjective and the constructed formal deformation will be versal.

In fact we claim there exists a sequence of deformations \ref{cart} with
the properties; 
\begin{list}{\textup{(\roman{temp})}}{\usecounter{temp}}
\item For each vertex $\{i\}$ there exists normal forms $$
\psi^{(n)}_{i}: \mathcal{R}_{\nu(\{i\})} \to R_{n}$$ 
lifting $\psi^{(n-1)}_{i}$ and such that the deformation $$
\Spec \Gamma(U_{i},
\mathcal{O}_{X_{n}})
\to \Spec R_{n} $$
of $Z_{\nu(\{i\})}$ is induced as in Proposition~\ref{nf} by
$\psi^{(n)}_{i}$.
\item Set $g_{ij}^{(n)} = \psi^{(n)}_{i}(t_{j}^{(1)})$ for all $i$ where
$\nu(\{i\}) = 6$, $\{j\} \in \link(\{i\})$ and let $\mathfrak{G}^{(n)}$ be
the set of these polynomials lifted to $P_{S}$. Then if
$\mathfrak{a}_{\mathfrak{G}^{(n)}}$ is as in \ref{aS} we have $R_{n+1}
\simeq P_{S}/(\mathfrak{a}_{\mathfrak{G}^{(n)}} +
\mathfrak{m}^{n+2})$.
\end{list} We start with the first-order case $n = 1$. For each $U_{i}$ we
exhibit the map $\mathcal{R}_{\nu(\{i\})}
\to R_{1}$ in Table~\ref{first-order}.  (With the convention when $\nu = 4$
that
 $t_{j}^{(0)}$ and $t_{k}^{(0)}$ (also $u_{i,j}$ and $u_{i,k}$) are the
same variable when $j$ and $k$ are opposite vertices in $\link(\{i\})$.)
Note that $g_{i,j}^{(1)} = \psi^{(1)}_{i}(t_{j}^{(1)}) = t_{i,j}$. 
\begin{table}
\centering
\begin{tabular}{l | l   l } 
Valency & $\mathcal{R}_{\nu(\{i\})}
\to R_{1}$ &\\
\hline $\nu(\{i\}) = 3$ & $t^{(-1)} \mapsto v_{i}$ & \\
 & $t_{j}^{(0)} \mapsto v_{i,j}, t_{j}^{(1)}
\mapsto t_{ij}$ & for each vertex $\{j\} \in
\link(\{i\})$ \\
 &$ t_{j}^{(2)} \mapsto v_{j,i}, t_{j}^{(3)} \mapsto v_{j}$ & if  $\{j\}
\in \link(\{i\})$ and $\nu(\{j\}) = 3$ \\
& $t_{j}^{(2)} \mapsto u_{j,i}$  & if  $\{j\} \in \link(\{i\})$ and
$\nu(\{j\}) = 4$\\
\hline $\nu(\{i\}) = 4$  & $t_{j}^{(0)} \mapsto u_{i,j}, t_{j}^{(1)}
\mapsto t_{ij}$ &  for each vertex $\{j\} \in
\link(\{i\})$\\
&$ t_{j}^{(2)} \mapsto v_{j,i}, t_{j}^{(3)} \mapsto v_{j}$ & if  $\{j\}
\in \link(\{i\})$ and $\nu(\{j\}) = 3$ \\
& $t_{j}^{(2)} \mapsto u_{j,i}$  & if  $\{j\} \in \link(\{i\})$ and
$\nu(\{j\}) = 4$\\
\hline $\nu(\{i\}) = 5,6$  & $t_{j}^{(1)}
\mapsto t_{ij}$ &  for each vertex $\{j\} \in
\link(\{i\})$\\
&$ t_{j}^{(2)} \mapsto v_{j,i}, t_{j}^{(3)} \mapsto v_{j}$ & if  $\{j\}
\in \link(\{i\})$ and $\nu(\{j\}) = 3$ \\
& $t_{j}^{(2)} \mapsto u_{j,i}$  & if  $\{j\} \in \link(\{i\})$ and
$\nu(\{j\}) = 4$\\
\hline
\end{tabular}
\caption{The first-order normal form for each $U_{i}$.}
\label{first-order}
\end{table}

Assume we have the deformations up to $R_{n}$. We must exhibit $X_{n+1}$
and the $\psi^{(n+1)}_{i}$ satisfying property (i) for the $R_{n+1}$
defined by property (ii). The ideal $\mathfrak{a}_{\mathfrak{G}^{(n)}}$
contains the images of the local obstruction equations \ref{obeq} for each
valency $6$ vertex. Thus each $\psi^{(n)}_{i}$ lifts to $\psi^{\prime}_{i}:
\mathcal{R}_{\nu(\{i\})} \to R_{n+1}$.  

Let $(U_{i},
\mathcal{O}^{\prime}_{i}) \to \Spec R_{n+1}$ be the induced normal form 
deformation of each chart. The difference between the deformations
$(U_{ij}, \mathcal{O}^{\prime}_{i})$ and $(U_{ij},
\mathcal{O}^{\prime}_{j})$ gives an element of $T^{1}_{U_{ij}}$. We know
that $H^{1}(\mathcal{T}_{S}^{1}) = 0$ (Theorem~\ref{T_man}), so we may
adjust these local deformations to make the difference 0. Explicitly we may
proceed as follows.

Recall that $U_{ij} = U_{i} \cap U_{j} = \emptyset$ if $\{i,j\}$ is 
not an edge. Assume that $\{i,j\}$ is an edge and that $\link(\{i,j\}) =
\{\{k\},
\{l\}\}$. In the local coordinates of $(U_{i},\mathcal{O}_{S})$ we may
write 
$$
\Gamma(U_{ij},\mathcal{O}_{S}) = k[y_{k}, y_{l}, y_{j},
y_{j}^{-1}]/(y_{k}y_{l})$$
where $y_{j}= x_{j}/x_{i}$ etc.. Thus we may represent the difference, i.e.
the element of $T^{1}_{U_{ij}}$, as 
$$
y_{k}y_{l} \mapsto \sum_{\alpha}a_{ij}^{\alpha}y_{j}^{\alpha}$$
with $a^{(\alpha)}_{ij} \in \mathfrak{m}^{n+1}/\mathfrak{m}^{n+2}$. 

If $\alpha
\ge 2$ change $\psi_{i}^{\prime}(t_{j}^{(\alpha)})$ to
$\psi_{i}^{\prime}(t_{j}^{(\alpha)}) -a_{ij}^{(\alpha)}$. If $\alpha
\le 0$ change $\psi_{j}^{\prime}(t_{i}^{(\alpha)})$ to
$\psi_{j}^{\prime}(t_{i}^{(\alpha)}) + a_{ij}^{(\alpha)}$. If $\alpha = 1$
we are free to adjust  $\psi_{i}^{\prime}(t_{j}^{(1})$ or
$\psi_{j}^{\prime}(t_{i}^{(1)})$ or both. Do this arbitrarily {\em unless}
one of the vertices, say $\{i\}$, has valency $6$ and the other not. In
this case adjust $\psi_{j}$ by adding $a_{ij}^{(1)}$ to the value of
$t_{i}^{(1)}$. 

Set $\psi^{(n+1)}_{i}$ to be the result after making these adjustments for
all edges $\{i,j\}$ and let $(U_{i},
\mathcal{O}^{(n+1)}_{i}) \to \Spec R_{n+1}$ be the new induced normal form 
deformation of each chart. The adjustments entail that for each $U_{ij}$ we
have isomorphisms $\phi_{ij} : \mathcal{O}^{(n+1)}_{i}|_{U_{ij}} \to
\mathcal{O}^{(n+1)}_{j}|_{U_{ij}}$. The next obstruction is in
$H^{2}(\mathcal{E}_{S}) = 0$ (Theorem~\ref{Da}). This means we may have to
adjust the $\phi_{ij}$, but not the $\mathcal{O}^{(n+1)}_{i}$, and
therefore not the normal form. We may now glue over these isomorphisms to
make $X_{n+1}$ with the wanted properties. 
\end{proof}

From Theorem~\ref{DXL_theorem} and the remark after it we get
\begin{corollary} \label{honest} There exists $\mathfrak{G}$ as in
Theorem~\ref{formal} 
and a local $k$-algebra $R$ with completion $\hat{R} =
\hat{R}_{\mathfrak{G}}$ such that $\Spec R$ is a versal base space for
$\Def^{a}_{S}$. In particular if all $g_{ij} = t_{ij}$ in $\mathfrak{G}$
then 
$ R =({P_{S}}/
\mathfrak{a}_{S})_{\mathfrak{m}} $.
\end{corollary}

An interesting set of examples comes about if we ask for all valencies for
vertices of $\mathcal{K}$ to equal $6$. This is known as a {\em degree $6$
regular triangulation}. If $n$ is the number vertices, then the the
$f$-vector must be $(n,3n,2n)$. In particular the Euler characteristic of
$\mathcal{K}$ is $0$ so $|\mathcal{K}|$ is a torus or a Klein bottle and
$S$ is either a degenerate abelian or bielliptic surface.

There are many such triangulations, see \cite{bk:equ} for a classification
for tori and  \cite{du:deg} for many examples. Certain such triangulations
where used to study degenerations of abelian surfaces in
\cite{gp:equ}.  We describe here just one series for the torus which
includes the vertex-minimal triangulation when $n=7$.
\begin{example}  \label{series} On $n$ vertices $\{0, \dots , n-1\}$ we
list the the $2n$ faces (all addition is done modulo $n$): 
$$
\{i, i + 2, i + 3\} \quad \{i, i + 1, i + 3\} \quad 0 \le i \le n-1 \, .$$
Note that $\link(\{i\})$ is the hexagon with vertices $\{i+2, i+3, i+1,
i-2, i-3, i-1\}$. This is the series $T_{n,1,2}$ in \cite{du:deg}.
\end{example}  It turns out that for such a triangulation we may choose all
$g_{i,j} = t_{i,j}$ in the description of the versal base space.
\begin{theorem} \label{torus} If $\mathcal{K}$ is a degree $6$ regular
triangulation of the torus or the Klein bottle and $ R =(k[t_{i,j} :
\{i,j\} \text{an edge in } \mathcal{K}]/
\mathfrak{a}_{S})_{\mathfrak{m}} $ then $\Spec R$ is a versal
base space for $\Def^{a}_{S}$. 
\end{theorem}
\begin{proof} We keep the notation from the proof of Theorem~\ref{formal}.
All $U_{i} \simeq Z_{6}$ and only the edges in $\mathcal{K}$ contribute to
$H^0(\mathcal{T}^1)$. Consider the equations 
in the normal form for a deformation of $Z_{6}$ with all $s_{j} = 0$; 
\begin{align*} &y_{j-1}y_{j+1}  +t^{(1)}_j y_{j}\quad j = 1, \dotsc ,6 \\
  &y_jy_{j+3}  -t^{(1)}_{j+1}t^{(1)}_{j+2} \quad  j = 1,2,3 \, .
\end{align*} 
For each $U_{i}$ we get a deformation in this normal form over the
completion $\hat{R}$ from the map $\psi_{i}:
\mathcal{R}_{6} \to \hat{R}$, $\psi_{i}(t_{j}^{(1)}) =  t_{ij}$ for each
each vertex $\{j\} \in
\link(\{i\})$. (Again we use the convention that the indices for $E_{6}$
are the indices of the vertices in $\link(\{i\})$ in cyclic order.) Let
$(U_{i}, \mathcal{O}_{i}) \to \Specf \hat{R}$ be the corresponding family.

We claim that we may construct a formal deformation 
\begin{equation*}
\begin{CD} X_{n} @>>> X_{n+1} \\
@VVV @VVV \\
\Spec\hat{R}/\mathfrak{m}^{n+1} @>>> \Spec\hat{R}/\mathfrak{m}^{n+2}
\end{CD}
\end{equation*} with $\Spec \Gamma(U_{i},
\mathcal{O}_{X_{n}}) = \mathcal{O}_{i}/\mathfrak{m}^{n+1} \mathcal{O}_{i} $
for all $n\ge 2$, i.e. at no level is it necessary to adjust the
$\psi_{i}$. Let $y_{j}^{(i)} = x_{j}/x_{i}$ be local coordinates for
$(U_{i},\mathcal{O}_{S})$. Assume that $\{i,j\}$ is an edge and that
$\link(\{i,j\}) = \{\{k\}, \{l\}\}$. We may write 
\begin{align*}
\Gamma(U_{ij},\mathcal{O}_{i}/\mathfrak{m}^{n+1} \mathcal{O}_{i})  & =
\hat{R}/\mathfrak{m}^{n+1}[y_{k}^{(i)}, y_{l}^{(i)}, y_{j}^{(i)},
(y_{j}^{(i)})^{-1}]/(y_{k}^{(i)}y_{l}^{(i)} + t_{ij}y_{j}^{(i)})\\
\Gamma(U_{ij},\mathcal{O}_{j}/\mathfrak{m}^{n+1} \mathcal{O}_{j}) & =
\hat{R}/\mathfrak{m}^{n+1}[y_{k}^{(j)}, y_{l}^{(j}, y_{i}^{(j)},
(y_{i}^{(j)})^{-1}]/(y_{k}^{(j)}y_{l}^{(j)} + t_{ij}y_{i}^{(j)}) \, .
\end{align*} Clearly $\phi_{ij}$ defined by $y_{k}^{(i)} \mapsto 
y_{k}^{(j)}/y_{i}^{(j)}$,  $y_{l}^{(i)} \mapsto 
y_{l}^{(j)}/y_{i}^{(j)}$ and  $y_{j}^{(i)} \mapsto 1/y_{i}^{(j)}$ is an
isomorphism. If $\{i,j,k\}$ is a face in $\mathcal{K}$ one checks that the
cocycle condition $\phi_{jk}\phi_{ij} = \phi_{ik}$ is satisfied on
$U_{ijk}$. 

Thus we have constructed a formal versal deformation over $\hat{R}$ and may
invoke Corollary~\ref{honest} to get the statement in the theorem.
\end{proof}

\begin{remark} The family constructed in the proof is only formal as one
can see by trying to make sense of the gluing isomorphisms over $U_{ijk}$
if $\{i,j,k\}$ is not a face. The line bundle $O_{S}(1)$ lifts trivially
over each power of $\mathfrak{m}$ so each $X_{r}$ is embedded via
 $y_{j}^{(i)} = x_{j}/x_{i}$ in $\mathbb{P}^{n}_{R/\mathfrak{m}^{r+1}}$,
but the equations defining $X_{r}$ are perturbed at each step. 
\end{remark}

\begin{example}  If $\mathcal{K}$ is one of the complexes in
Example~\ref{series} then $\mathfrak{a}_{S}$ is generated by the
minors of 
$$
\begin{bmatrix} t_{i,i+1} & t_{i,i+2} &  t_{i,i-3}\\
 t_{i,i-1} & t_{i,i-2} & t_{i,i+3}
\end{bmatrix} ,\quad i = 0, \dots n-1\, .$$
If $n=7$, i.e. we have the vertex-minimal triangulation of the torus then
the versal deformation has a very interesting structure involving a
6-dimensional reflexive polytope and a Calabi-Yau $3$-fold with Euler
number $6$. This will be studied in \cite{chr:ab}.\end{example}

\bibliographystyle{amsalpha}

\providecommand{\bysame}{\leavevmode\hbox to3em{\hrulefill}\thinspace}
\providecommand{\MR}{\relax\ifhmode\unskip\space\fi MR }
% \MRhref is called by the amsart/book/proc definition of \MR.
\providecommand{\MRhref}[2]{%
  \href{http://www.ams.org/mathscinet-getitem?mr=#1}{#2}
}
\providecommand{\href}[2]{#2}

\vspace{2ex}

{\small
\setbox0\hbox{FB Mathematik und Informatik, WE2}
\parbox{\wd0}{ Klaus Altmann\\
Institut f\"ur Mathematik\\
Freie Universit\"at Berlin\\
Arnimallee 3\\
D-14195 Berlin, Germany\\
email: altmann@math.fu-berlin.de}
\setbox1\hbox{email: christop@math.uio.no}\hfill
\parbox{\wd1}{ Jan Arthur Christophersen\\
Matematisk institutt\\
Postboks 1053 Blindern\\
University of Oslo\\
N-0316 Oslo, Norway\\
email: christop@math.uio.no}}

\end{document}